%% file: operator_Besov_eng.tex
\documentclass[11pt,reqno]{amsart}
\usepackage{amssymb,amscd,amsbsy,mathrsfs}
\input{classstartscr}

\usepackage{amsmath} 

\DeclareFontFamily{U}{mathx}{\hyphenchar\font45}
\DeclareFontShape{U}{mathx}{m}{n}{
  <5> <6> <7> <8> <9> <10>
        <10.95> <12> <14.4> <17.28> <20.74> <24.88>
        mathx10
}{}
\DeclareSymbolFont{mathx}{U}{mathx}{m}{n}
\DeclareFontSubstitution{U}{mathx}{m}{n}
\DeclareMathAccent{\widecheck}0{mathx}{"71}
\newcommand\OL{{\rm OL}}

\newcommand\cZ{\mathcal{Z}}
\newcommand\dg{\frak D}
\newcommand\ri{{\rm i}}

\newcommand\mB{\mathcal{B}}

\newcommand\mP{\mathcal{P}}

\newcommand\di{{\rm d}}

\newcommand{\Ast}{\mathbin{\widehat{\circledast}}}

\usepackage[T2A,T1]{fontenc}
\DeclareSymbolFont{cyrillic}{T2A}{cmr}{m}{it}
\def\makecyrsymbol#1#2{%
    \begingroup\edef\temp{\endgroup
        \noexpand\DeclareMathSymbol{\noexpand#1}
        {\noexpand\mathalpha}{cyrillic}%
        {\expandafter\expandafter\expandafter
            \calccyr\expandafter\meaning\csname T2A\string#2\endcsname\end}}%
    \temp}
\expandafter\def\expandafter\calccyr\string\char#1\end{#1}

\makecyrsymbol\Zhe\CYRZH
\makecyrsymbol\Be\CYRB
\makecyrsymbol\Shcha\CYRSHCH
\makecyrsymbol\Sha\CYRSH

\newcommand{\Bs}{\Be_{\be,1}^1}

\renewcommand{\C}{{\Bbb C}}
\renewcommand{\f}{{\varphi}}
\renewcommand{\b}{{\beta}}
\renewcommand{\L}{{\Lambda}}
\renewcommand{\l}{{\lambda}}

\renewcommand{\o}{{\omega}}
\renewcommand{\O}{{\Omega}}

\numberwithin{equation}{section}

\begin{document}

\noindent
\hspace*{5.5cm}To Oleg Vladimirrovich Besov
\newline
\hspace*{5.5cm}on the occasion of his 90th anniversary

\

\title{Besov spaces in operator theory}
\author{V.V. Peller}
\thanks{{The research on \S\:6-7  is supported by 
Russian Science Foundation [grant number 23-11-00153].}
}

\begin{abstract}
The survey is devoted to diverse applications of Besov classes in operator theory.
It is illustrated how Besov classes are used to describe Hankel operators of Schatten--von Neumann classes; various applications of this description are considered. Next, we discuss the role of Besov classes in norm estimates of polynomials of power bounded operators on Hilbert space and related estimates of Hankel matrices in tensor products of the spaces 
$\ell^1$ and
$\ell^\be$. 
An essential part of the survey is devoted to the role of Besov spaces in various problems of perturbation theory when studying the behavior of functions of a single operator or of a collection of operators under their perturbation. 
\end{abstract}

\maketitle

{\bf
\footnotesize
\tableofcontents
\normalsize
}

\

\setcounter{section}{0}
\section{\bf Introduction}
\setcounter{equation}{0}
\label{intr}

\

In this survey we discuss the important role played by the Besov classes in operator theory. We are going to demonstrate numerous situations when Besov classes are used in various issues of operator theory..

The Besov spaces were introduced in Besov's paper  \cite{Bes}. It turned out that they play a very noticeable role in embedding theorems, approximation theory, interpolation theorems and many other divisions of function theory, see \cite{BIN}, \cite{SNi}, \cite{Pee}, \cite{T}.

However, the importance of Besov spaces is not exhausted by the above.
Starting the 1980s, Besov classes have been playing a more and more important role in operator theory.
It is the role of Besov spaces in operator theory which is the principal subject of the survey.

We start with the Besov spaces $B_p^{1/p}(\T)$, $0<p<\be$, of functions on the unit circle $\T$, 
which appear while describing the Hankel operators of Schatten--von Neumann classes $\bS_p$.
We also consider certain applications and generalizations of these results such as commutators of multiplication operators and the Hilbert transform, Wiener--Hopf operators on a finite interval, the Ibragimov--Solev problem in the theory of stationary ransom processes.

Then we consider the problem of estimating the norms of polynomials of power bounded operators on Hilbert space. This problem is closely related to certain other Besov spaces of functions on the circle
$\T$ and with Hankel matrices in tensor products of the spaces  $\ell^1$ and 
$\ell^\be$. This leads in turn to a solution of Mazur's Problems 8 and 88 from the Scottish Book \cite{Ma}.

The final part of the survey is devoted to applications of Besov spaces in the theory of perturbations of operators on Hilbert space. We are going to consider various problems related to the behaviour of functions of operators or functions of collections of operators under perturbation. 

\

\pagebreak

\section{\bf Definitions of Besov classes}
\setcounter{equation}{0}
\label{Besov}

\

In this survey we are going to deal with Besov space on the unit circle (in other words, with Besov spaces of periodic function on the real line) and with Besov spaces on Euclidean spaces.
We start with the case of the unit circle. Moreover, we consider the case of the unit circle slightly in more detail
than the case of Euclidean space because Besov spaces on Euclidean spaces were discussed in many textbooks and monographs. We also observe that the definitions and properties of Besov spaces on the unit circle $\T$ 
can easily be generalized to the case of Besov spaces on the torus
 $\T^d$ for an arbitrary natural number $d$.

\subsection{Besov classes on the unit circle} Let $w$ be an infinitely differentiable function on  $\R$ such that
\bay
\label{w}
w\ge0,\quad\supp w\subset\left[1/2,2\right],\quad\mbox{and} \quad w(t)=1-w\left(\frac t2\right)\quad\mbox{for}\quad t\in[1,2].
\ey

We define the trigonometric polynomials $W_n$, $n\ge0$, by
\bay
\label{W_n}
W_n(\z)\df\sum_{j\in\Z}w\left(\frac{|j|}{2^n}\right)\z^j,\quad n\ge1,
\quad W_0(\z)\df\sum_{\{j:\,|j|\le1\}}\z^j,\quad\z\in\T.
\ey
If $f$ is a distribution (a generalized function) on $\T$, we put
$$
f_n=f*W_n,\quad n\ge0,
$$
and say that $f$ belongs to the {\it Besov class} $B_{p,q}^s(\T)$, $s\in\R$, 
$0<p,\,q\le\be$ if
\bay
\label{Bperf}
\big\{2^{ns}\|f_n\|_{L^p}\big\}_{n\ge0}\in\ell^q.
\ey

Recall that a distribution on $\T$ is a continuous linear functional on the space of infinitely differentiable functions on the circle. The class $L^1(\T)$ can naturally be imbedded in the space of distributions: a function $g$ in $L^1(\T)$ determines the linear functional
$f\mapsto\int_\T fg\,d\m$, where $\m$ is normalized Lebesgue measure on $\T$.

In the case when $p,q\ge1$, $B_{p,q}^s(\T)$ is a Banach space. It can be endowed with various equivalent norms. In particular, we can put
\bay
\label{pqsnorma}
\|f\|_{B_{p,q}^s(\T)}=\Big\|\{2^{ns}\|f_n\|_{L^p}\big\}_{n\ge0}\Big\|_{\ell^p}.
\ey 
If $\min\{p,q\}<1$, the expression \rf{pqsnorma} is a quasinorm. Note that there are other natural norms (quasinorms) on Besov classes.

If $p=q$, it is common to use the notation $B_p^s(\T)\df B_{p,p}^s(\T)$.

The dual space with the Besov space $B_{p,q}^s(\T)$ in the case when
 $p,q\ge1$ and $q<\be$, can be described in the following way:
$$
\Big(B_{p,q}^s(\T)\Big)^*=B_{p',q'}^{-s}(\T),\qquad p'=p/(1-p),\quad q'=q/(1-q),
$$
with respect to the bilinear form 
$$
(f,g)=\sum_{j\in\Z}\hat f(j)\hat g(j);
$$
here $f$ ranges over the dense subset of trigonometric polynomials in $B_{p,q}^s(\T)$ while $g$ is an arbitrary function  (distribution) in $B_{p',q'}^{-s}(\T)$. Recall that if $g$ is a distribution on $\T$, the Fourier coefficient $\hat g(j)$ is the value of the functional $g$ at the function $\bar z^j$.

Let us proceed now to the description of the Besov spaces $B_{p,q}^s(\T)$ for $s>0$ in terms of finite differences. Suppose for simplicity that $p\ge1$ and $q\ge1$. We define the difference operator $\D_\t$, $\t\in\T$, on the set of functions on $\T$ by
$$
(\D_\t f)(\z)\df f(\t\z)-f(-z),\quad\z\in\T,
$$
while the difference operator $\D^n_\t$ is defined as the $n$th power of $\D_\t$.

Suppose that $s>0$, $1\le p\le\be$ and $1\le q\le\be$. Then the following description of the spaces 
$B_{p,q}^s(\T)$ holds:
$$
B_{p,q}^s(\T)=\left\{f\in L^p(\T):~\int_\T
\frac{\|\D^n_\t f\|_{L^p}^q}{|1-\t|^{1+sq}}\,\di \m(\t)<\be
\right\},
\quad q<\be;
$$
$$
B_{p,\be}^s(\T)=
\left\{
f\in L^p(\T):~\sup_{\t\ne1}\frac{\|\D^n_\t f\|_{L^p}}{|1-\t|^s}<\be
\right\}.
$$

In the case $p=q=\be$ and $0<s<1$, we obtain the H\"older class $\L_s$:
$$
f\in\L_s=B_{\be,\be}^s(\T)\quad\Longleftrightarrow\quad|f(\z_1)-f(\z_2)|\le\const|\z_1-\z_2|^s,\quad\z_1,\z_2\in\T.
$$
The classes $\L_s\df B_\be^s(\T)=B_{\be,\be}^s(\T)$ for $s>0$ form the so-called H\"older--Zygmund
scale of spaces. In particular, $\L_1$ is the Zygmund space:
$$
f\in\L_1\quad\Longleftrightarrow\quad|f(\z\t)-2f(\z)+f(\z\bar\t)|\le\const|1-\t|,\quad\z,\t\in\T.
$$
We denote by $\l_s$ the {\it closure of the trigonometric polynomials in} $\L_s$
and we say that $\l_s$ is the {\it separable H\"older--Zygmund space}.

We proceed now to the description of Besov spaces in terms of harmonic extension to the unit disc. 
Let $s\in\R$
and let $m$ be a nonnegative integer such that $m>s$.
Then a distribution $f$ on the circle $\T$ belongs to the space
$B_{p,q}^s(\T)$  if and only if
\bay
\label{q<be}
\int_0^{1}(1-r)^{(m-s)q-1}\left\|\frac{\partial^m}{\partial r^m}\Big((\mP f)(r\z)\Big)\right\|_{L^p(\T)}^q\,\di r
<\be,\quad q<\be,
\ey
\bay
\label{q=be}
\sup_{r\in[0,1)}(1-r)^{m-s}\left\|\frac{\partial^m}{\partial r^m}\Big((\mP f)(r\z)\Big)\right\|_{L^p(\T)}<\be,\quad q=\be,
\ey
where $\mP f$ denotes the Poisson integral of the distribution $f$.

Consider the Riesz projections $\pp_+$ and $\pp_-$ defined on the class of all distributions on  $\T$ by
\bay
\label{Riesz}
\pp_+f=\sum_{j\ge0}\hat f(j)z^j,\quad\pp_-f=\sum_{j<0}\hat f(j)z^j.
\ey

It follows easily from the definition of Besov classes in terms of the convolutions with $W_n$ that
$$
\pp_+B_{p,q}^s(\T)\subset B_{p,q}^s(\T)\quad\mbox{and}\quad
\pp_-B_{p,q}^s(\T)\subset B_{p,q}^s(\T),\qquad p,q>0,\quad s\in\R.
$$
We are going to use the notation $\big(B_{p,q}^s\big)_+$ for $\pp_+B_{p,q}^s(\T)$.
Functions $f$ in $\big(B_{p,q}^s\big)_+$ \lb(distributions for $s\le0$) can be identified with the power series
$$
\sum_{j\ge0}\hat f(j)z^j
$$
which makes the classes $\big(B_{p,q}^s\big)_+$ spaces of analytic functions in the disc. Such classes can be described in terms of conditions \rf{q<be} and \rf{q=be}; herewith the partial derivatives of order $m$ can be replaced with the derivatives of order $m$ in complex variable. In particular, for $p=q<\be$ and $s<0$ the Besov classes  
$\big(B_p^s\big)_+\df\big(B_{p,p}^s\big)_+$ become  {\it weighted Bergman classes}, that is a function  $f$ analytic in the unit disc $\dd$ belongs to the Besov class
$\big(B_p^s\big)_+$ if and only if
$$
\int_\dd(1-|\z|)^{-sp-1}| f(\z)|^p\,\di \m_2(\z)<\be,\quad p<\be.
$$
The important Bloch space ${\frak B}$ of functions $f$ analytic in $\dd$ such that
$$
\sup_{\z\in\dd}|f'(\z)|(1-|\z|)<\be,
$$
is also a Besov class, namely $\big(B_\be^0\big)_+$, of analytic functions.

Let us mention here the book \cite{Pe}, which contains an introduction to Besov spaces on the unit circle.

\subsection{Besov classes on Euclidean spaces} Let $d$ be a positive integer and let $w$ be an infinitely smooth function on $\R$ satisfying requirements \rf{w}.
We define the functions $W_n$, $n\in\Z$, on $\R^d$ by 
$$
\big(\F W_n\big)(x)=w\left(\frac{\|x\|}{2^n}\right),\quad n\in\Z, \quad x=(x_1,\cdots,x_d),
\quad\|x\|\df\Big(\sum_{j=1}^dx_j^2\Big)^{1/2},
$$
where $\F$ is {\it Fourier transform} defined on $L^1\big(\R^d\big)$ by
$$
\big(\F f\big)(t)=\!\int\limits_{\R^d} f(x)e^{-{\rm i}(x,t)}\,dx,\!\quad 
x=(x_1,\cdots,x_d),
\quad t=(t_1,\cdots,t_d), \!\quad(x,t)\df \sum_{j=1}^dx_jt_j.
$$
Clearly,
$$
\sum_{n\in\Z}(\F W_n)(t)=1,\quad t\in\R^d\setminus\{0\}.
$$
Note that the utilization of the same notation $W_n$ as in previous subsections will not lead to a confusion.

With each tempered distribution $f$ in ${\mathscr S}^\prime\big(\R^d\big)$ we associate the sequence 
$\{f_n\}_{n\in\Z}$,
\bay
\label{fn}
f_n\df f*W_n.
\ey
The formal series $\sum_{n\in\Z}f_n$, being a Littlewood--Payley type expansion of $f$, does not necessarily converges to $f$. 
We start with the definition of the (homogeneous) Besov class $\dot B^s_{p,q}\big(\R^d\big)$,
$s\in\R$, $0<p,\,q\le\be$, as the space of distributions
$f$ of class ${\mathscr S}^\prime(\R^n)$ (i.e., tempered distributions) such that
\bay
\label{Wn}
\{2^{ns}\|f_n\|_{L^p}\}_{n\in\Z}\in\ell^q(\Z),\quad
\|f\|_{B^s_{p,q}}\df\big\|\{2^{ns}\|f_n\|_{L^p}\}_{n\in\Z}\big\|_{\ell^q(\Z)}.
\ey
In accordance with this definition, $\dot B^s_{p,q}(\R^d)$ contains all polynomials; next, $\|f\|_{B^s_{p,q}}=0$ for an arbitrary polynomial $f$. Moreover, the distribution $f$ is determined by the sequence $\{f_n\}_{n\in\Z}$
uniquely up to a polynomial. It is easy to see that the series 
$\sum_{n\ge0}f_n$ converges in\footnote{Here and in what follows we assume that the space 
${\mathscr S}^\prime(\R^d)$ is equipped with the weak topology
$\s\big({\mathscr S}^\prime(\R^d),{\mathscr S}(\R^d)\big)$.} ${\mathscr S}^\prime(\R^d)$. 
However, the series $\sum_{n<0}f_n$ can diverge in general. Nevertheless, it can be shown that the series
\bay
\label{ryad}
\sum_{n<0}\frac{\partial^r f_n}{\partial x_1^{r_1}\cdots\partial x_d^{r_d}}\qquad \mbox{for}\quad r_j\ge0,\quad
1\le j\le d,\quad\sum_{j=1}^dr_j=r,
\ey
converge uniformly in $\R^d$, if $r\in\Z_+$ and
$r>s-d/p$. Note that for $q\le1$, the series in \rf{ryad}
converge uniformly under the weaker condition $r\ge s-d/p$.

Now we can define the modified (homogeneous) Besov space $B^s_{p,q}\big(\R^d\big)$. We say that a distribution $f$
belongs to the Besov class $B^s_{p,q}(\R^d)$ if condition \rf{Wn} is satisfied and
$$
\frac{\partial^r f}{\partial x_1^{r_1}\cdots\partial x_d^{r_d}}
=\sum_{n\in\Z}\frac{\partial^r f_n}{\partial x_1^{r_1}\cdots\partial x_d^{r_d}}\quad
\mbox{for}\quad 
r_j\ge0,\quad\mbox
1\le j\le d,\quad\sum_{j=1}^dr_j=r,
$$
in the space ${\mathscr S}^\prime\big(\R^d\big)$, where $r$ is the minimal nonnegative integer, for which
$r>s-d/p$ ($r\ge s-d/p$ if $q\le1$). Now the function $f$ is uniquely determined by the sequence $\{f_n\}_{n\in\Z}$ up to a polynomial of degree less than $r$. Moreover, a polynomial $g$ belongs to 
$B^s_{p,q}\big(\R^d\big)$ if and only if $\deg g<r$. 

In the case $p=q$ we use the notation $B_p^s(\R^d)$ for $B_{p,p}^s(\R^d)$.

Consider now the scale $\L_\a(\R^d)$, $\a>0$,  of {\it H\"older--Zygmund classes}. They can be defined as $\L_\a(\R^d)\df B_\be^\a(\R^d)$.

The Besov classes admit many other descriptions. Let us state the description in terms of finite differences.
For $h$ in $\R^d$, we define the difference operator $\D_h$ by
$(\D_hf)(x)=f(x+h)-f(x)$, $x\in\R^d$.

Let $s>0$,  $m\in\Z$ and $m-1\le s<m$.  
Let $p,q\in[1,+\be]$.
The Besov class $B_{p,q}^s\big(\R^d\big)$ can be defined as the space of
functions $f$ in $L^1_{\rm loc}\big(\R^d\big)$ such that
$$
\int_{\R^d}|h|^{-d-sq}\|\D^m_h f\|_{L^p}^q\,\di h<\be,~\; q<\be;\qquad
\sup_{h\not=0}\frac{\|\D^m_h f\|_{L^p}}{|h|^s}<\be,~\; q=\be.
$$
However, with this definition the Besov classes can contain polynomials of degrees higher than in the case of the definition in terms of convolutions with the functions $W_n$.

The space $B_{pq}^s(\R^d)$ can be defined in terms of Poisson integral.
Let $P_d(x,t)$ be the Poisson kernel that solves the Dirichlet problem for the half-space 
$\R^{d+1}_+\lb\df\{(x,t): x\in\R^d, t>0\}$, i.e., $P_d(x,t)=c_dt(|x|^2+t^2)^{-\frac{d+1}2}$,
$c_d=\pi^{-\frac{d+1}2}\G(\frac{d+1}2)$. With each function 
$f\in L^1\big(\R^d,(\|x\|+1)^{-(d+1)}\,dx\big)$
we can associate the Poisson integral $\mP f$,
$$
(\mP f )(x,t)=\int_{\R^d}P_d(x-y,t)f(y)\,\di y.
$$
Then for each positive integer $m$, the following equality holds
$$
\frac{\partial^m (\mP f)}{\partial^m t}(x,t)=\int_{\R^d}\frac{\partial^m P_d(x-y,t)}{\partial^m t}f(y)\,\di y.
$$
Note that the second integral makes sense for all 
$f$ in $L^1\big(\R^d,(\|x\|+1)^{-(d+m+1)}\,\di x\big)$.
This allows one to define $\frac{\partial^m}{\partial t^m}\mP$ for functions $f$ of class 
$L^1\big(\R^d,(\|x\|+1)^{-(d+m+1)}\,\di x\big)$.

Let $m\in\Z$, $s>0$, \mbox{$m-1\le s<m$}, $1\le p,\:q\le+\be$. Then we can define the space
$B_{pq}^s(\R^d)$ as the space of functions 
$f\in L^1\big(\R^d,(\|x\|+1)^{-(d+m+1)}\,dx\big)$ such that
$$
\left(\int_0^{\be}t^{(m-s)q-1}\left\|\Big(\frac{\partial^m}{\partial t^m}\mP f\Big)(\cdot,t)\right\|_{L^p(\R^d)}^q\,\di t\right)^{\frac1q}
<+\be,\quad q<+\be,
$$
$$
\sup_{t>0}t^{m-s}\left\|\Big(\frac{\partial^m}{\partial t^m}\mP f\Big)(\cdot,t)\right\|_{L^p(\R^d)}<+\be,\quad
q=+\be.
$$
Under this definition Besov classes can also contain polynomials of higher degrees than in the case of the definition in terms of convolutions with the functions $W_n$. Note also that the definition in terms of Poisson integral under certain conditions and provisions can work even in the case when $p<1$  or $q<1$.

\subsection{Inhomogeneous Besov classes on Euclidean spaces.} As a rule, it is homogeneous Besov spaces that appear in operator theory. However, in the recent  papers \cite{APne} and \cite{APdi} inhomogeneous Besov spaces also appear.

We define the function $W^{[0]}$ by
$$
\big(\F W^{[0]}\big)(t)=1-\sum_{n\ge1}(\F W_n)(t),\quad t\in\R^d. 
$$
It is easy to see that
$$
(\F W^{[0]}\big)(t)=1,\quad\mbox{if}\quad\|t\|\le1\quad\mbox{and}\quad
\supp\F W^{[0]}\subset\big\{t\in\R^d:~\|t\|\le2\big\}.
$$
Let $s>0$ and $p,q\in[1,\be]$.
We say that a tempered distribution $f$ on $\R^d$ belongs to the {\it inhomogeneous Besov space} $\Be^s_{p,q}(\R^d)$ if
\bay
\label{fnol'}
f^{[0]}\df f*W^{[0]}\in L^p(\R^2)\quad\mbox{and}\quad
\big\{2^{ns}\|f_n\|_{L^p(\R^d)}\big\}_{n\ge1}\in\ell^q,
\ey
where the functions $f_n$ are defined by \rf{fn}.

Note that this notation for inhomogeneous Besov classes was adopted in  \cite{APne} and \cite{APdi}, and its is somewhat different from the one commonly used. Nevertheless, it is convenient.

It is well known and it is easy to see that $\Be^s_{p,q}(\R^d)=B^s_{p,q}(\R^d)\bigcap L^p(\R^d)$. 

We refer the reader to the books \cite{Pee}, \cite{BIN}, \cite{SNi} and \cite{T} for more detailed information on Besov spaces.

\

\section{\bf An introduction to Schatten--von Neumann classes}
\setcounter{equation}{0}
\label{SchatNeu}

\

In this section we give a brief introduction in singular values of Hilbert space operators and the Schatten--von Neumann classes. We refer the reader to \cite{GK} and \cite{BS} for the proofs of the facts given below and for more detailed information.

For a bounded linear operator $T$ on Hilbert space, its  {\it singular values} $s_j(T)$, $j\ge0$, are defined by 
\bay
\label{sj}
s_j(T)\df\inf\big\{\|T-R\|:~\rank R\le j\big\}.
\ey
In the case when $T$ is a compact operator on Hilbert space, the sequence 
$\{s_j(T)\}_{j\ge0}$ coincides with the sequence of eigenvalues (counted with multiplicities) of the operator 
$(T^*T)^{1/2}$ arranged in the nonincreasing order.

{\it The Schatten--von Neumann class} $\bS_p$, $0<p<\be$, consists, by definition, of the operators
$T$, for which
\bay
\label{normaSp}
\|T\|_{\bS_p}\df\Big(\sum_{j\ge0}\big(s_j(T)\big)^p\Big)^{1/p}<\be.
\ey
For $p\ge1$, this is a normed ideal of operators on Hilbert space with norm \rf{normaSp}.

For $p<1$, the functional $\|\cdot\|_{\bS_p}$ is not a norm, it is a $p${\it-norm}, i.e.,
the following inequality holds:
\bay
\label{pnorma}
\|T+R\|_{\bS_p}^p\le\|T\|_{\bS_p}^p+\|R\|_{\bS_p}^p,\quad T,~R\in\bS_p.
\ey

The class $\bS_1$ is called {\it trace class}. If $T$ is a trace class operator on Hilbert space $\h$, its {\it trace} 
$\trace T$ is defined by
$$
\trace T\df\sum_{j\ge0}(Te_j,e_j),
$$
where $\{e_j\}_{j\ge0}$ is an orthonormal basis in $\h$. The right-hand side does not depend on the choice of a basis.

The class $\bS_2$ is called the {\it Hilbert--Schmidt class}. It forms a Hilbert space with inner product
$$
(T,R)_{\bS_2}\df\trace(TR^*).
$$

For $p\in(1,\be)$, the dual space $(\bS_p)^*$ can isometrically be identified with the space 
$\bS_{p'}$, $1/p+1/p'=1$, with respect to the bilinear form
$$
\langle T,R\rangle\df\trace(TR).
$$
Moreover, for $R\in\bS_{p'}$, the norm of the linear functionsl $T\mapsto\trace(TR)$ on $\bS_p$ coincides with the norm of $R$ in $\bS_{p'}$.

The space dual to $\bS_1$ can be identified with the space of all bounded linear operators with respect to the same bilinear form while the space dual to the space of compact operators  can be identified with $\bS_1$.


\

\section{\bf Hankel operators and Besov classes}
\setcounter{equation}{0}
\label{Hankel}

\

\subsection{Classical Hankel operators and Hankel matrices}
Besov classes were used to study Hankel operators in \cite{PeSb1}.
Hankel operators on the Hardy class\footnote{We refer the reader to \cite{Gar}, \cite{Ku} and \cite{Nik} for definitions and basic properties of the Hardy classes $H^p$ as well as the classes ${\rm BMO}$ and ${\rm VMO}$.} 
$H^2$ of functions analytic in $\dd$
can be defined in the following way. Let $\f$ be a function of class $L^2=L^2(\T)$. On the set of analytic polynomials $\cp_+$ the Hankel operator $H_\f$ from the Hardy space $H^2$ 
to the space $H^2_-\df L^2\ominus H^2$ is defined by
$$
H_\f\df\pp_-\f f,\quad f\in\cp_+,
$$
(recall that the Riesz projections are defined by \rf{Riesz}\,). Such operators are not necessarily bounded. The boundedness criterion for Hankel operators is given by the Nehari theorem:
the operator $H_\f$ is bounded if and only if there exists a function $\psi$ in $L^\be$
such that $H_\f=H_\psi$. The last equality is equivalent to the following one:
$$
\hat\f(j)=\hat\psi(j),\quad j<0.
$$
Herewith
$$
\|H_\f\|=\dist_{L^\be}(\psi,H^\be).
$$
It follows from the theorems of Nehari and Ch. Fefferman that for $\f\in L^2$, the Hankel operator $H_\f$ is a bounded operator from $H^2$ to $H^2_-$ (to be more precise, it extends to a bounded operator from $H^2$ to $H^2_-$) if and only if the function $\pp_-\f$ belongs to the class ${\rm BMO}$ of functions on $\T$  of bounded mean oscillation. Herewith the norm of $H_\f$ is equivalent to the norm
of $\pp_-\f$ in the space ${\rm BMO}$. 
 
A compactness criterion for Hankel operators is given by a the Hartman theorem: if $\f\in L^\be$, then  $H_\f$ is compact if and only if $\f\in H^\be+C$, i.e., $\f$ can be represented as
$\f=g+h$, where $g\in H^\be$ and $g$ is a continuous function on $\T$. Together with a theorem by Sarason this leads to the following criterion: the operator $H_\f$ is compact if and only if the function $\pp_-\f$ belongs to the space ${\rm VMO}$ of functions on $\T$ of vanishing mean oscillation.

Detailed information on Hankel operators and their applications can be found in the book  \cite{Pe}, see also the book \cite{Nik}.

The Hankel operator $H_\f:H^2\to H^2_-$ in the orthonormal bases $\{z^j\}_{j\ge0}$ of the space $H^2$
and $\{\bar z^j\}_{j>0}$ of the space $H^2_-$ has infinite matrix
$$
\left(
\begin{matrix}
\hat\f(-1)&\hat\f(-2)&\hat\f(-3)&\cdots\\
\hat\f(-2)&\hat\f(-3)&\cdots&\cdots\\
\hat\f(-3)&\vdots&\ddots\\
\vdots
\end{matrix}
\right).
$$
Such matrices are called {\it Hankel matrices}. We associate with a sequence of complex numbers 
$\{\a_j\}_{j\ge0}$ in $\ell^2$ the Hankel matrix
\bay
\label{Ga}
\G_\a=\left(
\begin{matrix}
\a_0&\a_1&\a_2&\cdots\\
\a_1&\a_2&\cdots&\cdots\\
\a_2&\vdots&\ddots\\
\vdots
\end{matrix}
\right)
\ey
and we are going to identify such matrices with operators on the space $\ell^2$ with respect to the natural basis in $\ell^2$. Clearly, the study of metric properties of Hankel operators $H_\f$ is equivalent to the study of metric properties of Hankel matrices $\G_\a$.
In particular, $\G_\a$ induces a bounded operator on $\ell^2$ if and only if
\bay
\label{stryad}
\a(z)\df\sum_{j\ge0}\a_jz^j\in{\rm BMOA}\df\pp_+{\rm BMO}={\rm BMO\cap H^2},
\ey
and in the case when the last condition is satisfied, the following equality holds
$$
\|\G_\a\|=\|\a\|_{\rm BMOA}\df\inf\{\|\psi\|_{L^\be}:~\psi\in L^\be\quad\mbox{and}\quad\hat\psi(j)=\a_j\quad\mbox{for}\quad j\ge0\}.
$$


On the other hand, the operator $\G_\a$ is compact if and only if
$$
\a(z)\df\sum_{j\ge0}\a_jz^j\in{\rm VMOA}\df\pp_+{\rm VMO}={\rm VMO\cap H^2}.
$$

We proceed now to conditions under which Hankel operators belong to the Schatten--von Neumann classes
$\bS_p$. The following theorem describe the Hankel operators $H_\f$ that belong to the Schatten--von Neumann classes
$\bS_p$.

\begin{thm}
\label{HfSp}
Let $0<p<\be$ and let $\f\in L^2$. Then the Hankel operatos $H_\f$ belongs to the Schatten--von Neumann class $\bS_p$ if and only if the function $\pp_-\f$ belongs to the Besov class
$B_p^{1/p}(\T)$.
\end{thm}

Theorem \ref{HfSp} can be reformulated in the following way:
$$
\G_\a\in\bS_p\qquad\Longleftrightarrow\qquad\sum_{j\ge0}\a_jz^j\in\big(B_p^{1/p}\big)_+,\qquad0<p<\be.
$$

Applying the description of Besov classes \rf{q<be} to Besov classes of analytic functions, we obtain the following description of the operators $\G_\a$ of class $\bS_p$. 

\begin{thm}
Let $\a$ be given by the power series {\em\rf{stryad}}. Let $0<p<\be$ and let
$m$ be an integer such that $m>1/p$. Then
$$
\G_\a\in\bS_p\qquad\Longleftrightarrow\qquad\int_\dd|\a^{(m)}(\z)|^p(1-|\z|)^{mp-2}\,\di \m_2(\z)<\be.
$$
\end{thm}

In particular, taking $p=1$, we obtain the following description of trace class operators $\G_\a$:

\begin{cor}
$$
\G_\a\in\bS_1\qquad\Longleftrightarrow\qquad\int_\dd|\a''(\z)|\,\di \m_2(\z)<\be.
$$
\end{cor}

Note that the description of the trace class Hankel operators was obtained in  \cite{PeSb1}. This solved an old well known problem; in particular, it was mentioned by M. Rosenblum \cite{Ro}, M.G. Krein \cite{Kr} and F. Holland (see \cite{ABB}).

The case $1<p<\be$ was also covered in \cite{PeSb1}.  In the case $0<p<1$ the description of Hankel operators of class $\bS_p$ was obtained in \cite{PeSb2}. Note that other approaches in the case 
$0<p<1$ we offered by A.A. Pekarskii \cite{Pek} and S. Semmes \cite{Se}.

The framework of this survey do not allow to give here a complete proof of Theorem  \ref{HfSp}. We restrict ourselves to a proof of the sufficiency of the condition $\pp_-\f\in B_p^{1/p}(\T)$ for $0<p<1$ and also we discuss how it is possible to obtain the necessity of this condition in the case $p=1$. 

\medskip

{\bf Proof of the sufficiency of the condition $\bs{\pp_{-}\f\in B_p^{1/p}(\T)}$ for $\bs{p\le1}$.} Let 
$\psi$ be a function analytic in the disc $\dd$ and let $\psi_n=\psi*W_n$, where the analytic polynomials
$W_n$ are defined by \rf{W_n}. Using inequalty \rf{pnorma} (which is obviously also valid for $p=1$), we see that it follows from \rf{Bperf} and \rf{pqsnorma} that it suffices to establish the following assertion:

{\it Let $\psi$ be an analytic polynomial of degree $m-1$. Then
\bay
\label{degpsi}
\|\G_\psi\|_{\bS_p}\le2^{1/p-1}m^{1/p}\|\psi\|_{L^p(\T)}.
\ey}

Now, for an analytic polynomial $\psi$, $\deg\psi=m-1$, we define the operators $A_j$, \lb$0\le j\le2m-1$,
of rank 1 as follows. Consider the points $\z_j=e^{2\pi\ri j/(2m)}$, $0\le j\le2m-1$,
on the circle $\T$.
Let $f_j$ and $g_j$, $0\le j\le2m-1$, be sequences in $\ell^2$
defined by
$$
f_j(k)=\left\{\begin{array}{ll}
\ov{\psi(\z_j)}\z_j^k,&0\le k<m,\\
0,&k\ge m;
\end{array}\right.
\qquad
g_j(n)=\left\{\begin{array}{ll}
\bar\z_j^n,&0\le k<m,\\
0,&k\ge m.
\end{array}\right.
$$
We define now the rank one operator $A_j$, $0\le j\le2m-1$, by
$$
A_jx=(x,f_j)g_j,\quad x\in\ell^2.
$$
Then the following equality holds
\bay
\label{predGpsi}
\G_\psi=\frac1{2m}\sum_{j=0}^{2m-1}A_j.
\ey
To prove this equality, we have to verify that the matrix entries of the operators on the left and on the right are the same. Clearly, in the case when $k\ge m$ or $n\ge m$, we have
$$
(A_je_k,e_n)=0\quad\mbox{and}\quad(\G_\psi e_k,e_n)=\hat\psi(n+k)=0.
$$
Here $\{e_j\}_{j\ge0}$ is the standard orthonormal basis in $\ell^2$. Next, since $\psi$ is an analytic polynomial and $\deg\psi<m$, the following identity holds
$$
\frac1{2m}\sum_{j=0}^{2m-1}\psi(\z_j)\bar\z_j^d=\hat\psi(d)\quad\mbox{for}\quad 0\le d\le 2m-1.
$$ 
Thus, for
$n<m$ and $k<m$ we have $(A_je_k,e_n)=\psi(\z_j)\bar\z_j^{n+k}$, and so
$$
\frac1{2m}\sum_{j=0}^{2m-1}(A_je_k,e_n)=\frac1{2m}\sum_{j=0}^{2m-1}\psi(\z_j)\bar\z_j^{n+k}=\hat\psi(n+k).
$$
Now, representation \rf{predGpsi} together with inequality \rf{pnorma} allows us to conclude that
$$
\|\G_\psi\|^p_{\bS_p}\le\frac1{(2m)^p}\sum_{j=0}^{2m-1}\|A_j\|^p_{\bS_p}
=\frac1{(2m)^p}\sum_{j=0}^{2m-1}\|f_j\|^p\|g_j\|^p
=\frac1{(2m)^p}\sum_{j=0}^{2m-1}(|\psi(\z_j)|m)^p,
$$
that is
\bay
\label{diskLp}
\|\G_\psi\|_{\bS_p}\le 2^{1/p-1}m^{1/p}\frac1{2m}\sum_{j=0}^{2m-1}|\psi(\z_j)|^p.
\ey
To obtain \rf{degpsi}, it suffices to use the well known inequality
$$
\frac1{2m}\sum_{j=0}^{2m-1}|\psi(\z_j)|^p\le\|\psi\|_{L^p}^p.
$$
Incidentally, we can simply apply inequality \rf{diskLp} for the function $\psi_\t$, $\t\in\T$, $\psi_\t(\z)\df\psi(\t\z)$, observe that $\|\G_{\psi_\t}\|_{\bS_p}=\|\G_\psi\|_{\bS_p}$, $\t\in\T$, and integrate it in the variable $\t$ with respect to Lebesgue measure $\m$ on $\T$. $\bl$

\medskip

For $p=1$, there are several approaches to the proof of the necessity of the condition $\a\in\big(B_1^1\big)_+$. One of them is selecting special bounded operators $B$ an $\ell^2$ and calculating the  traces $\trace\G_\a B$ (see \cite{PeSb1} and \cite{Pe}). Another one is constructing a certain special bounded projection ${\mathscr Q}$ from $\bS_1$ onto the subspace of trace class Hankel operators (see \cite{Pe}, Ch. 6, \S\:5).

Note here that the averaging projection ${\mathcal P}$ onto the Hankel matrices defined by 
\bay
\label{proekt}
({\mathcal P}T)_{jk}\df\frac1{j+k+1}\sum_{m=0}^{j+k}(Te_m,e_{j+k-m}),\quad j,~k\ge0,
\ey
is not a bounded operator on $\bS_1$. Indeed, it is easy to observe that as 
$T$ ranges over the set of projections of rank 1, the Hankel matrices $\cp T$ fill all matrices of the form
$\G_\a$, where $\a$ is an arbitrary analytic in $\dd$ function with derivative in the Hardy class $H^1$ and so 
$$
{\mathcal P}\bS_1=\big\{\G_\f:~\f'\in H^1\big\}.
$$
Nevertheless, 
$\{f:~f'\in H^1\}\not\subset\big(B_1^1\big)_+$.
 Indeed, it is well known that the lacunary Fourier coefficients $\hat g(2^k)$ of the functions $g$ in $H^1$ fill the whole space $\ell^2$ while  the lacunary Fourier coefficients of the functions with derivatives in $(B_1^1\big)_+$
fill the space $\ell^1$. 

We refer the reader to \cite{PeSb2} and \cite{Pe} for the proof of the necessity of the condition $\pp_-\f\in B_p^{1/p}(\T)$, $0<p<1$, in Theorem \ref{HfSp}. The proof is rather difficult and is based on explicit computations of the $\bS_p$-quasinorms of convolution operators on cyclic groups.

In the case $1<p<\be$, the proof of Theorem \ref{HfSp} given in the paper \cite{PeSb1} was obtained with the help of the Marcinkiewicz interpolation theorem. Moreover, the reasoning given in the proof implies that the averaging projection ${\mathcal P}$ onto the Hankel matrices is bounded on $\bS_p$ for $1<p<\be$. We refer the reader to the book  \cite{Pe}, in which all such questions are discussed in detail.

\subsection{Generalized Hankel matrices} Let $\a$ and $\b$ be real numbers and let $\f$ be a function analytic in the disc $\dd$. Consider {\it generalized Hankel matrices} of the form
$$
\G_\f^{\a,\b}=\big\{(1+j)^\a(1+k)^\b\big\}_{j,k\ge0}.
$$
The case $\a=\b=0$ corresponds to the classical Hankel matrices. We are going to consider such matrices as operators on $\ell^2$. 

The following theorem describes the bounded and compact operators of the form $\G_\f^{\a,\b}$ for $\a\ge0$ and $\b\ge0$.

\begin{thm}
Let $\f$ be a function analytic in $\dd$. The following statements hold:

{\em(i)} if $\a>0$ and $\b>0$, then $\G_\f^{\a,\b}$ is a matrix of a bounded operator on $\ell^2$ if and only if 
$\f$ belongs to the H\"older class $\L_{\a+\b}$;

{\em(ii)} if $\a>0$ and $\b>0$, then $\G_\f^{\a,\b}$ is a matrix of a compact operator on $\ell^2$ if and only if 
$\f$ belongs to the separable H\"older class $\l_{\a+\b}$;

{\em(iii)} if $\a>0$, then $\G_\f^{\a,0}$ is a matrix of a bounded operator on $\ell^2$ if and only if the function 
$$
I_{-\a}\f\df\sum_{n\ge0}(1+n)^\a\hat\f(n)z^n
$$
belongs to ${\rm BMO}$;

{\em(iv)} if $\a>0$, then $\G_\f^{\a,0}$ is a matrix of a compact operator on $\ell^2$ if and only if the function  $I_{-\a}\f$ belongs to ${\rm VMO}$.
\end{thm}

The following theorem describes the operators $\G_\f^{\a,\b}$ in the Schatten--von Neumann class $\bS_p$.

\begin{thm}
Let $\f$ be a function analytic in $\dd$ and let $0<p<\be$. Suppose that $\a$ and $\b$ are real numbers such that 
$$
\min\{\a,\b\}>\max\Big\{-\frac12,-\frac1p\Big\}.
$$
Then $\G_\f^{\a,\b}\in\bS_p$ if and only if $\f$ belongs to the Besov class
$\big(B_p^{1/p+\a+\b}\big)_+$.
\end{thm}

We refer the reader to the book \cite{Pe} for proofs of the above assertions and for more detailed discussions.

\

\section{\bf Certain applications of the criterion of the membership of Hankel operators in classes   
$\boldsymbol{\bS_p}$}
\setcounter{equation}{0}
\label{pioH}

\

In this section we consider certain applications of the description of the Hankel operators of class $\bS_p$ in terms of Besov classes. The reader can learn applications in more detail in the book \cite{Pe}. 

\subsection{A characterization of Besov classes in terms of rational approximation}

In this subsection we discuss how the Besov classes  $\big(B_p^{1/p}\big)(\T)$ can be described
in terms of rational approximation in the norm of the space ${\rm BMO}$. 

First of all, we need the oldest theorem on Hankel matrices, that is the Cronecker theorem \cite{Kro} on
finite rank Hankel matrices, which implies that a Hankel operator $H_\f:H^2\to H^2_-$ has rank $k$ if and only if  $\pp_-\f$ is a rational function of degree $k$.\footnote{The degree $\deg p/q$ of the rational function $p/q$ in the case when $p$ and $q$ are polynomials that have no nontrivial common divisor is defined as 
$\max\{\deg p,\deg q\}$.}

We are also going to use the profound Adamyan--Arov--Krein theorem \cite{AAK}, which asserts that
in the case when $T$ is a Hankel operator, it suffices to consider on the right-hand side of
\rf{sj} exclusively Hankel operators $R$ of rank at most $j$, i.e.,
$$
s_j(H_\f)=\inf\big\{\|H_\f-H_g\|:~g\in L^\be,~\rank H_g\le j\big\}.
$$

A combination of the Kronecker  theorem with the Adamyan--Arov--Krein theorem together with the Nehari and the Fefferman theorems mentioned above allows us to obtain the following description of the Besov classes
$\big(B_p^{1/p}\big)(\T)$ in terms of rational approximation in the norm of  ${\rm BMO}$.

Let $\f\in L^\be$ and $j\in\Z_+$. Consider the deviations $\rho_n(\f)$ of the function $\f$ from the set
${\mathcal R}_n$ of rational functions of degree at most $n$ with poles off $\T$ in the norm of ${\rm BMO}$:
$$
\rho_n(\f)\df\dist_{\rm BMO}(\f,{\mathcal R}_n).
$$

\begin{thm}
\label{ratsap}
Let $\f$ be a function $\T$ of class ${\rm BMO}$ and let $0<p<\be$. Then $\f\in\big(B_p^{1/p}\big)(\T)$
if and only if
$$
\sum_{n\ge0}(\rho_n(\f))^p<\be.
$$
\end{thm}

Theorem \ref{ratsap} was obtained in \cite{PeSb1} and \cite{PeSb2}. We refer the reader to \cite{Pe},
Ch. 6, \S\:6 where these issues are discussed in detail and other related results are given.

\subsection{Commutators of multiplication operators and the Riesz projection} For a function $\f$ of class $L^2$, the {\it multiplication operator} $M_\f$ is defined on the dense subset $\cp$ of  trigonometric polynomials by $M_\f f\df\f f$, $f\in\cp$. Consider now the {\it commutator} ${\mathscr C}_\f$ of the operator $M_\f$ and Riesz projection $\pp_+$ defined by
$$
{\mathscr C}_\f=[M_\f,\pp_+]f=\f\pp_+f-\pp_+\f f,\quad f\in\cp.
$$
Such commutators play an important role in harmonic analysis. It is well known that 
$$
(\pp_+f)(\z)={\rm v.p.}\int_\T\frac{f(\t)}{1-\bar\t\z}\,\di \m(\t)+\frac12f(\z),\quad f\in\cp,\quad\z\in\T.
$$
Here the notation ${\rm v.p.}$ means the integral in the sense of principal value. It is easy to deduce from the above expression the following identity
\bay
\label{komCf}
{\mathscr C}_\f f={\rm v.p.}\int_\T\frac{\f(\z)-\f(\t)}{1-\bar\t\z}f(\t)\,\di \m(\t),\quad f\in\cp.
\ey

Let us express now the commutator ${\mathscr C}_\f$ in terms of Hankel operators. It is convenient to use the notation $f_+\df\pp_+f$ and $f_-\df\pp_-f$.

\begin{thm}
Let $\f\in L^2$ and $f\in\cp$. Then
$$
{\mathscr C}_\f f=H_{\f_-} f_+-\big(H_{\ov{\f_+}}\big)^*f_-.
$$
\end{thm}

\Pf We have
\begin{align*}
{\mathscr C}_\f f&=\f\pp_+f-\pp_+\f f=\f_+f_++\f_-f_+-\f_+f_+-\pp_+\f_-f_+-\pp_+\f_+f_-\\[.2cm]
&=\pp_-\f_-f_+-\pp_+\f_+f_-=H_{\f_-} f_+-\big(H_{\ov{\f_+}}\big)^*f_-.\quad\bl
\end{align*}

\medskip

Since the Hankel operator $H_{\f_-}$ acts from $H^2$ to $H^2_-$ and the operator 
$\big(H_{\ov{\f_+}}\big)^*$ acts from $H^2_-$ to $H^2$, it is easy to see that the commutator ${\mathscr C}_\f$ is bounded (compact) if and only if both Hankel operators $H_{\f_-}$ and $H_{\ov{\f_+}}$ are bounded (compact). This implies the well known results in harmonic analysis: the boundedness of the commutator
${\mathscr C}_\f$ is equivalent to the membership of $\f$ in the class ${\rm BMO}$  while the compactness of 
${ \mathscr C}_\f$ s equivalent to the membership of $\f$ in the class ${\rm VMO}$.

In the same way we can consider the problem of the membership of the commutator ${\mathscr C}_\f$ 
in the Schatten--von Neumann classes. Thus, the following criterion holds.

\begin{thm}
\label{Cfi}
Let $0<p<\be$. Then the commutator ${\mathscr C}_\f$ defined by {\em\rf{komCf}} belongs to the class 
$\bS_p$ if and only if  $\f$ belongs to the Besov class $\big(B_p^{1/p}\big)(\T)$.
\end{thm}

We refer the reader to the book \cite{Pe} for more detailed information.

\subsection{Hankel operators as integral operators on $\bs{L^2[0,\be)}$; commutators of multi\-plication operators and the Hilbert transform} For a function $k$ of class $L^1(\R_+)$ the integral Hankel operator $\bs\varGamma_k$ on the space $L^2(\R_+)$ is defined by
$$
(\bs\varGamma_kf)(t)=\int_0^\be k(s+t)f(s)\,\di s.
$$
The operators $\bs\varGamma_k$ are integral analogues of the matrix operators $\G_\a$ defined by \rf{Ga}. It turns out, however that the operators $\bs\varGamma_k$ are not only integral analogs of Hankel operators  $\G_\a$ but also have Hankel matrices in the basis of Laguerre functions. We are going to discuss in this subsection how to define the operators $\bs\varGamma_k$ in the situation considerably more general than in the case $k\in L^1(\R_+)$, and we are also going to discuss the conditions of boundedness, compactness,  membership in the Schatten--von Neumann ideals for such operators. 

To state the results, we need the notions of a distribution and a tempered distribution.

For an open subset $\O$ of the real line $\R$, we consider the space
${\mathcal D}(\O)$ of infinitely differentiable functions on $\O$ with compact support and we equip this space with  
the natural topology of inductive limit. The {\it space ${\mathcal D}'(\O)$ of distributions} on $\O$ is defined as the space of continuous linear functionals on ${\mathcal D}(\O)$. We refer the reader to the L. Scwartz book   \cite{Sch} for an introduction to the theory of distributions.

To define the space of tempered distributions, we consider the Scwartz space  ${\mathcal S}$ that consists of infinitely differentiable functions $f$ on $\R$ such that
$$
\sup_{t\in\R}\big|f^{(m)}(t)\big|(1+|t|)^n<\be\qquad\mbox{for all}\qquad m\quad\mbox{and}\quad n\quad\mbox{in}\quad\Z_+.
$$
We can define now  {\it the space ${\mathcal S}'$ of tempered distributions} as the space of continuous linear functionals on the space ${\mathcal S}$ endowed with the natural topology, see \cite{Sch}.

It is well known that the Fourier transform $\F$ and the inverse inverse Fourier transform  $\F^*$ defined on the space  $L^1(\R)$ by
$$
(\F f)(s)=\int_\R e^{-2\pi{\ri}ts}f(t)\,\di t\quad\mbox{and}\quad(\F^* f)(t)=\int_\R e^{2\pi{\ri}ts}f(s)\,\di s,
$$
map the Schwartz space ${\mathcal S}$ onto itself which makes it possible to extend them by duality to 
linear maps of the space of tempered distributions ${\mathcal S}'$ onto itself.

The following result holds:

\begin{thm}
Let $k\in{\mathcal D}'(0,\be)$. The following statements are equivalent:

{\em1)} $\bs\varGamma_k$ is a bounded operator on $L^2[0,\be)$;

{\em2)} there exists a function $\vk$ of class $L^\be(\R)$ such that $(\F\vk)\big|(0,\be)=k$;

{\em3)} there exists a distribution $r$ in ${\mathcal S}'$ such that $r\big|(0,\be)=k$,
$\supp r\subset[0,\be)$ and $\F^*r\in{\rm BMO}(\R)$.
\end{thm}

We refer the reader to the book \cite{Pe}, Ch. 1, \S\:8, where such issues are discussed in detail and, in particular, a compactness criterion for the integral operators $\bs\varGamma_k$ is given.

Consider now the operator ${\mathcal U}$ defined on the space $L^2(\T)$ by
$$
({\mathcal U}f)(t)=\frac1{\sqrt{\pi}}\frac{(f\circ\o^{-1})(t)}{t+\ri},\quad t\in\R,
$$
where
$$
\o(\z)\df\ri\frac{1+\z}{1-\z},\quad\z\in\dd.
$$
It is well known that ${\mathcal U}$ is a unitary operator from $L^2(\T)$ onto $L^2(\R)$; moreover,
It maps the Hardy class $H^2$ on the disc onto the Hardy class $H^2(\C_+)$ on the half-plane.

The following assertion holds:

\begin{thm}
The bounded operators $\bs\varGamma_k$ have Hankel matrices in the orthonormal basis  
$\{\F{\mathcal U}z^n\}_{n\ge0}$ in $L^2(\R_+)$. Moreover, the following identity holds:
$$
\big(\F{\mathcal U}z^n\big)(t)=-2\sqrt{\pi}\ri L_n(4\pi t)e^{-2\pi t},\quad t>0,
$$
where the $L_n$ are the Laguerre polynomials defined by 
$$
L_n(s)=\frac1{n!}e^s\Big(\frac{d}{ds}\Big)^n(e^{-s}s^n).
$$
\end{thm}

The proofs of the above statements can be found in \cite{Pe}, Ch. 1, \S\:8.

Let us now state a a necessary and sufficient condition for an integral operator $\bs\varGamma_k$ to belong to the Schatten--von Neumann classes $\bS_p$.

\begin{thm}
\label{intGS_p}
Let $0<p<\be$ and let $k$ be a distribution on $(0,\be)$. Then $\bs\varGamma_k\in\bS_p$
if and only if $k=(\F\psi)\big|(0,\be)$ for a function $\psi$ in the Besov class 
$B_p^{1/p}(\R)$.
\end{thm}

We refer the reader for details to the book \cite{Pe}, Ch. 6, \S\:7. 

Theorem \ref{intGS_p} allows us to obtain an analog of Theorem \ref{Cfi} in the case of commutators of multiplication and the Hilbert transform. For a function $\f$ of class ${\rm BMO}(\R)$, we consider the operator 
$\bs{\mathcal C}_f$ defined by
$$
(\bs{\mathcal C}_\f f)(y)=\mbox{v.p.}\int_\R\frac{\f(x)-\f(y)}{x-y}f(x)\,\di x,\quad f\in L^2(\R).
$$
It is well known that the operator $\bs{\mathcal C}_\f$ is bounded on $L^2(\R)$ if and only if $\f\in{\rm BMO}(\R)$.

\begin{thm}
Let $0<p<\be$ and let $\f$ be a function of class ${\rm BMO}(\R)$. Then the operator $\bs{\mathcal C}_\f$ belongs to the Schatten--von Neumann class $\bS_p$ if and only if $\f$ belongs to the Besov class $B_p^{1/p}(\R)$.
\end{thm}

\subsection{Wiener--Hopf operators on a finite interval} 
Let $\s>0$. We define the Wiener--Hopf operators on the space $L^2[-\s,\s]$ in the following way. Let $k$ be a locally integrable function on  $(-2\s,2\s)$ (or even a distribution on $(-2\s,2\s)$).
Consider the operator $W_{\s,k}$ defined on the dense subset  in $L^2[-\s,\s]$ of infinitely smooth functions $f$ with compact support in $(-\s,\s)$ by
$$
(W_{\s,k}f)(x)=\int_{-\s}^\s k(x-y)f(y)\,\di y,\quad x\in(-\s,\s).
$$
Such operators are called {\it Wiener--Hopf operators on a finite interval}. In this subsection we discuss conditions, under which such operators belong to $\bS_p$ for $0<p<\be$.

Together with the operators $W_{\s,k}$ we are going to consider so-called {\it truncated Hankel operators} $\bs{\varGamma}_{\s,k}$ on $L^2[-\s,\s]$ defined by
$$
(\bs{\varGamma}_{\s,k}f)(x)=\int_{-\s}^\s k(x+y)f(y)\,\di y,\quad x\in(-\s,\s).
$$
It is easy to see that $\bs{\varGamma}_{\s,k}=W_{\s,k}U$ where $U$ is the unitary operator on $L^2[-\s,\s]$ defined by
$$
(Uf)(x)=f(-x),\quad x\in(-\s,\s).
$$

To state the description of the operators $W_{\s,k}$ and $\bs{\varGamma}_{\s,k}$ of class $\bS_p$, we define the following functions:
$$
\nu_j(x)=w\left(\frac{x+2\s}{2^{j+1}\s}\right),\quad j<0;\qquad\nu_j(x)=\nu_{-j}(-x),\quad j>0;\qquad
\nu_0(x)=1-\sum_{j\ne0}\nu_j(x),
$$
where $w$ is an infinitely smooth function satisfying \rf{w}.

The following theorem was proved in \cite{PeWH} (note that in \cite{Roc} the same result was established by different methods in the case when $p\ge1$).

\begin{thm}
Let $\s>0$, $0<p<\be$, and let $k$ be a distribution on $(-2\s,2\s)$. 
The following statements are equivalent:

{\em(i)} $W_{\s,k}\in\bS_p$;

{\em(ii)}$\bs{\varGamma}_{\s,k}\in\bS_p$;

{\em(iii)}
$$
\sum_{j\in\Z}2^{-|j|}\|\F(\nu_jk)\|_{L^p}^p<\be
$$

{\em(iv)} there exist functions $\psi_1$ and $\psi_2$ in the Besov class $B_p^{1/p}(\R)$ such that
$$
(\F\psi_1)(x)=k(x-2\s),\quad x\in(0,4\s),
$$
and
$$
(\F\psi_2)(x)=k(2\s-x),\quad x\in(0,4\s).
$$
\end{thm}

\

In conclusion we note that there are many generalizations of the results in subsections  4 and 5, in which Besov classes appear in the description certain operators that belong to $\bS_p$. Some of those results are listed in the references of the book \cite{Pe}. However, besides them there are many other publications that appeared after the publication of \cite{Pe}.

\subsection{The I.A. Ibragimov--V.N. Solev problem on stationary random processes} The Ibragimov--Solev problem posed in the note \cite{IS2} concerns regularity properties of stationary processes. The problem is 
to find out when the product of the projections onto the Past and onto the Future belongs to the Schatten--von Neumann class $\bS_p$, $1\le p<\be$.

Without entering into details of the theory of stationary processes, we note that we understand by a (wide sense) stationary process a sequence of random variables 
 $\{X_n\}_{n\in\Z}$ such that the mathematical expectation ${\Bbb E}X_jX_k$ depends only on the difference $j-k$. Let 
$$
c_n\df{\Bbb E}X_nX_0.
$$
By the F. Riesz--Herglotz theorem (see \cite{RSN}, \S\:53),  there exists a finite positive Borel measure $\mu$ on the unit circle $\T$ such that $c_n=\int_\T z^n\,d\mu$. The measure $\mu$ is called the {\it spectral measure of the process}. Then the closed linear span
$\clos\spn\{X_n:~n\in\Z\}$ can be isometrically identified with the weighted space $L^2(\mu)$ via the map 
$$
X_j\mapsto z^j,\quad j\in\Z.
$$
Herewith the past of the process is identified with the subspace 
$$
H^2_-(\mu)\df\clos\spn\{z^j:~j<0\},
$$
while the future starting the moment $n$ is identified with the subspace
$$
z^nH^2(\mu)\df\clos\spn\{z^j:~j\ge n\}.
$$

We are going to consider {\it regular stationary processes}, i.e., such processes, for which 
$$
\bigcap_{n\ge0}z^nH^2(\mu)=\{\0\}.
$$
It is well known that a stationary process is regular if and only if its spectral measure
 $\mu$ is absolutely continuous with respect to Lebesgue measure and its density
 $w$ (it is called the {\it the spectral density of the process}) satisfies the condition $\log w\in L^1(\T)$. 

Suppose that we have a regular stationary process with spectral density $w$. 
Consider now the orthogonal projections ${\frak P}_0$ and ${\frak P}^n$ onto the past and the future starting with  moment $n$. The Ibragimov--Solev problem was to describe the spectral densities of regular stationary processes, for which the operators  $\cp_0\cp^n$ belong to the Schatten--von Neumann class $\bS_p$. The following theorem gives a solution to this problem, see \cite{PeSb1} and \cite{PeKh}.

\begin{thm}
Let $1\le p<\be$ and let $w$ be the spectral density of a regular stationary process. The following are equivalent:

{\em(i)} ${\frak P}_0{\frak P}^n\in\bS_p$ for all $n\in\Z_+$;

{\em(ii)} ${\frak P}_0{\frak P}^0\in\bS_p$;

{\em(iii)} $w$ admits a representation 
$$
w=|P|^2e^\f,
$$
where $P$ is a polynomial with zeros on $\T$ and $\f$ is a real-valued function in the Besov class 
$B_P^{1/p}(\T)$.
\end{thm}

Note that for $p=2$ this assertion was proved earlier in the paper by Ibragimov and Solev \cite{IS1} and the methods of \cite{IS1} do not work for $p\ne2$.

\

\section{\bf Power bounded operators and the tensor products 
$\bs{\ell^\be\widehat{\otimes}\ell^\be}$ 
and 
$\bs{\ell^\1\widecheck{\otimes}\ell^\1}$; 
Masur's problems from the Scottish Book}
\setcounter{equation}{0}
\label{tens}

\

This section is devoted to finding estimates of polynomials of power bounded operators on Hilbert space. 
It turns out that the resulting norms are closely related to the tensor products $\ell^\be\widehat{\otimes}\ell^\be$ adn $\ell^\1\widecheck{\otimes}\ell^\1$. Herewith the problem of the description of Hankel matrices of class
$\ell^\1\widecheck{\otimes}\ell^\1$ appears and we obtain such a description in terms of a Besov class. We also study properties of the averaging projection onto the space of Hankel matrices in such norms. In conclusion we will see that such results lead to the solution of Mazur's problem 8 and 88 from the Scottish book.

\subsection{Estimates of operator polynomials}
We  start this subsection with the problem of estimating functions of power bounded operators on Hilbert space.
We are going to discuss certain such estimates. In particular, we will see that the norm of a polynomial of a power bounded operator can be estimated in terms of the norm of the polynomial in the Besov class
  $\big(B_{\be,1}^0\big)_+$.

Thus, let $T$ be a power bounded operator on Hilbert space, i.e.,
\bay
\label{Tnc}
\|T^n\|\le c,\quad n\ge0.
\ey
Here $c>1$ (the case $c=1$ corresponds to the case when $T$ is a contraction, i.e., $\|T\|\le1$). It is well known that for $c>1$, the operator $T$ does not have to be {\it polynomially bounded}, i.e., 
for such operators $T$ the estimate
$$
\|\f(T)\|\le\const\max_{|\z|\le1}|\f(\z)|
$$ 
for analytic polynomials $\f$ does not hold in general
(this was established in \cite{Leb}, see also \cite{Pe82} and the papers quoted there).

The paper \cite{Pe82} was devoted to the problem of finding most optimal estimates for the norms
$\|\f(T)\|$ under condition \rf{Tnc}. In other words, the problem is to estimate the norm 
$$
|||\f|||_c\df\sup\big\{\|\f(T)\|:~T~\mbox{ satisfies the condition~ \rf{Tnc}}\big\}
$$
in an optimal way.

To state the results, we introduce certain tensor products. We define {\it the projective tensor product } $\ell^\be\widehat\otimes\ell^\be$ as the space of matrices 
$\{\g_{jk}\}_{j,k\ge0}$ that admit a representation
\bay
\label{prlbe}
\g_{jk}=\sum_n f_n(j)g_n(k),
\ey
where $\{f_n(j)\}_{j\ge0}\in\ell^\be$, $\{g_n(k)\}_{k\ge0}\in\ell^\be$ and
\bay
\label{tennorm}
\sum_n\|f_n\|_{\ell^\be}\|g_n\|_{\ell^\be}<\be.
\ey
By the norm $\big\|\{\g_{jk}\}_{j,k\ge0}\big\|_{\ell^\be\widehat\otimes\ell^\be}$ of the matrix 
$\{\g_{jk}\}_{j,k\ge0}$ in $\ell^\be\widehat\otimes\ell^\be$ we understand the infimum of the left-hand side of inequality 
\rf{tennorm} over all representations of this matrix in the form \rf{prlbe}. Equipped with this norm, the space $\ell^\be\widehat\otimes\ell^\be$ becomes a Banach space.

Consider now the tensor algebra $V^2$ (the Varopoulos algebra) as a weak completion of the space 
$\ell^\be\widehat\otimes\ell^\be$ in the following sense. We define the projections $P_m$, $m\in\Z_+$, on the space of all infinite matrices by
$$
(P_m\g)_{jk}=\left\{\begin{array}{ll}\g_{jk},&j\le m,~k\le m\\[.2cm]
0,&\mbox{otherwise}.
\end{array}\right.
$$
We say that an infinite matrix $\g=\{\g_{jk}\}_{j,k\ge0}$ {\it belongs to the space $V^2$} if
$$
\sum_{m\in\Z_+}\|P_m\g\|_{\ell^\be\widehat\otimes\ell^\be}<\be.
$$
Herewith $\ell^\be\widehat\otimes\ell^\be$ is a closed proper subspace of $V^2$. It is easy to observe that the identity matrix belongs to $V^2$ but does not belong to $\ell^\be\widehat\otimes\ell^\be$.

It is easy to see that if $A=\{a_{jk}\}_{j,k\ge0}\in V^2$, then $A$ is a {\it Schur multiplier} of the space
${\mathcal B}(\ell^2)$ of matrices of bounded linear operators on $\ell^2$, i.e.,
$$
B=\{b_{jk}\}_{j,k\ge0}\in{\mathcal B}(\ell^2)\quad\Longrightarrow\quad
A\star B\df\{a_{jk}b_{jk}\}_{j,k\ge0}\in{\mathcal B}(\ell^2).
$$
It turns out that the converse also holds, see \cite{Be}.

Let us now define the {\em injective tensor product} $\ell^1\widecheck\otimes\ell^1$ of two copies of $\ell^1$ as the completion of the algebraic tensor product 
$\ell^1\otimes\ell^1$ in the norm
$$
\big\|\{\g_{jk}\}_{j,k\ge0}\big\|_{\ell^1\widecheck\otimes\ell^1}=
\sup\left\{\left|\sum_{j,k\ge0}\g_{jk}f_jg_k\right|:
~\|\{f_j\}_{j\ge0}\|_{\ell^\be}\le1,~\|\{g_k\}_{k\ge0}\|_{\ell^\be}\le1\right\}.
$$
It well known that the dual space $(\ell^1\widecheck\otimes\ell^1)^*$ can be identified with the space  $V^2$ with respect to the bilinear form
$$
\big(\{\g_{jk}\}_{j,k\ge0},\{\b_{jk}\}_{j,k\ge0}\big)=\sum_{j,k\ge0}\g_{jk}\b_{jk},\quad
\{\g_{jk}\}_{j,k\ge0}\in\ell^1\widecheck\otimes\ell^1,\quad\{\b_{jk}\}_{j,k\ge0}\in V^2.
$$
Note also that the space $\ell^1\widecheck\otimes\ell^1$ can naturally be identified with the space of bounded linear operators from the space $c_0$ of sequences with zero limit to the space $\ell^1$.

Let us now define the space ${\mathscr L}$ of functions $f$ analytic in the disc $\dd$, for which there exists a matrix $\{\g_{jk}\}_{j,k\ge0}$ in $\ell^1\widecheck\otimes\ell^1$ such that
\bay
\label{norL}
\hat f(n)=\sum_{j+k=n}\g_{jk},\quad n\ge0.
\ey
The norm $\|f\|_{{\mathscr L}}$ of a function $f$ in the space ${\mathscr L}$ is defined as the infimum of the norms $\|\{\g_{jk}\}_{j,k\ge0}\|_{\ell^1\widecheck\otimes\ell^1}$
over all matrices $\{\g_{jk}\}_{j,k\ge0}$ in $\ell^1\widecheck\otimes\ell^1$ that satisfy \rf{norL}. It is easy to show that 
$$
\|\f\|_{{\mathscr L}}=\sup\big\{|\langle\f,\psi\rangle|:~\|\G_\psi\|_{V^2}\le1\big\},
$$
that is, the supremum is taken over all analytic in $\dd$ functions $\psi$, for which the Hankel matrix $\G_\psi$ belongs to the unit ball of the space $V^2$. Here
$$
\langle\f,\psi\rangle\df\sum_{n\ge0}\hat\f(n)\hat\psi(n).
$$ 

To give the estimate obtained in \cite{Pe82} of the norms of polynomials of power bounded operators 
in terms of the norm of the polynomial in ${\mathscr L}$, we recall {\it Grothendieck's inequality} \cite{Gr}.
Let $\{\g_{jk}\}_{j,k\ge0}\in\ell^1\widecheck\otimes\ell^1$, $\{x_j\}_{j\ge0}$ and $\{y_k\}_{k\ge0}$ be vectors in the unit ball of a Hilbert space. Then
\bay
\label{neGr}
\left|\sum_{j,k\ge0}\g_{jk}(x_j,y_k)\right|\le k_{\rm G}\|\{\g_{jk}\}_{j,k\ge0}\|_{\ell^1\widecheck\otimes\ell^1},
\ey
where $k_{\rm G}$ is the so-called Grothendieck constant, i.e., the least number, for which the inequality holds.

The following result was obtained in \cite{Pe82}.

\begin{thm}
\label{otsL}
Let $T$ be an operator on Hilbert space that satisfies {\em\rf{Tnc}}.
Then
\bay
\label{cherezL}
\|\f(T)\|\le c^2k_{\rm G}\|\f\|_{{\mathscr L}}
\ey
for every analytic polynomial $\f$.
\end{thm}

\Pf Let $x$ and $y$ be vectors in the unit ball of the Hilbert space and let $\{\g_{jk}\}_{j,k\ge0}$ be a complex matrix whose entries except for finitely many are equal to zero and such that 
$$
\sum_{j+k=n}\g_{jk}=\hat\f(n),\quad n\ge0.
$$
We have
$$
(\f(T)x,y)=\sum_{n\ge0}\hat\f(n)(T^nx,y)=\sum_{j,k\ge0}\hat\f(j+k)(T^jT^kx,y)
=\sum_{j,k\ge0}\g_{jk}(T^kx,(T^*)^jy).
$$
Since $\|T^kx\|\le c$ and $\|(T^*)^jy\|\le c$, we obtain from Grothendieck's inequality \rf{neGr} 
$$
|(\f(T)x,y)|\le c^2k_{\rm G}\|\{\g_{jk}\}_{j,k\ge0}\|_{\ell^1\widecheck\otimes\ell^1}
$$
which proves the theorem. $\bl$

Let us proceed now to more tangible estimates of the norms of operator polynomials.

\medskip

{\bf Definition.} We say that a function $\psi$ analytic in $\dd$ is a {\it multiplier of the Hardy class} $H^1$ if
$$
f\in H^1\quad\Longrightarrow\quad\psi*f\df\sum_{n\ge0}\hat f(n)\hat\psi(n)z^n\in H^1.
$$

\medskip

We denote by ${\frak M}H^1$ the class of multipliers of $H^1$. By the norm of a function $\psi$ in 
${\frak M}H^1$ we mean the norm of the operator $f\mapsto\psi*f$ on $H^1$.
It is easy to see that $\psi\in{\frak M}H^1$ if and only if $\psi$ is a multiplier of the space ${\rm BMOA}$, see \rf{stryad}.

Theorem \ref{otsL} implies easily the following estimate:

\begin{cor}
Let $T$ be a power bounded operator on Hilbert space.
Then
\bay
\label{MH1}
\|\f(T)\|\le\const\sup\big\{|\langle\f,\psi\rangle|:~\|\psi\|_{{\frak M}H^1}\le1\big\}
\ey
for every analytic polynomial $\f$.
\end{cor}

Let us define now the space ${\rm VMOA}\Ast H^1$ of functions analytic in the disc $\dd$ by
$$
{\rm VMOA}\Ast H^1=\left\{\sum_mf_m*g_m:~f_m\in{\rm VMOA},~g_m\in H^1,~
\sum_m\|f_m\|_{\rm VMOA}\|g_m\|_{H^1}<\be
\right\}.
$$
The norm in ${\rm VMOA}\Ast H^1$ is defined by
$$
\|\f\|_{{\rm VMOA}\Ast H^1}\df\inf\left\{
\sum_m\|f_m\|_{{\rm VMOA}}\|g_m\|_{H^1}:~\sum_mf_m\ast g_m=\f
\right\}.
$$
It is easy to see that
\bay
\label{VMOH1}
\big({\rm VMOA}\Ast H^1\big)^*={\frak M}H^1
\ey
with respect to the natural pairing.

Now we are in a position to state the result of \cite{Pe82} on estimates of polynomials of power bounded operators in terms of the norms of the polynomials in ${\rm VMOA}\Ast H^1$.

\begin{thm}
\label{otsAst}
Let $T$ be a power bounded operator on Hilbert space. Then
$$
\|\f(T)\|\le\const\|\f\|_{\tiny{{\rm VMOA}\Ast H^1}}
$$
for every analytic polynomial $\f$.
\end{thm}

Theorem \ref{otsAst} follows easily from \rf{MH1} and \rf{VMOH1}.

Finally, we are able to estimate the norms of polynomials of power bounded operators in terms of the norms of
the polynomials in the Besov space $\big(B^0_{\be,1}\big)_+$.

\begin{thm}
\label{otsBes}
Let $T$ be a power bounded operator on Hilbert space. Then
$$
\|\f(T)\|\le\const\|\f\|_{B^0_{\be,1}}
$$
for every analytic polynomial $\f$.
\end{thm}

\Pf By Theorem \ref{otsAst}, it suffices to prove that
$$
\big(B^0_{\be,1}\big)_+\subset{\rm VMOA}\Ast H^1
$$
and the embedding is continuous.

\medskip

Let $\f\in\big(B^0_{\be,1}\big)_+$. It is easy to see that
$$
\f=\sum_{n\ge0}\f*W_n=\f*V_0+\sum_{n\ge1}(\f*W_n)*(W_{n-1}+W_n+W_{n+1}).
$$
Thus,
$$
\|\f\|_{{\rm VMOA}\Ast H^1}\le|\hat\f(0)|+\sum_{n\ge1}\|\f*W_n\|_{L^\be}\|W_{n-1}+W_n+W_{n+1}\|_{L^1}.
$$
The result follows now from the well known fact that $\|W_n\|_{L^1}\le\const$. $\bl$

\medskip

A natural question arises of whether the above estimates for the norms of polynomials of power bounded operators are optimal.
Suppose that $\|\cdot\|_\flat$ is a norm on the set of analytic polynomials and we have the estimate
$$
\|\f(T)\|\le\const\|\f\|_\flat
$$
for an arbitrary operator $T$ satisfying \rf{Tnc}. Also, the constant on the right-hand side of \rf{Tnc}  may depend on the number $c$ in \rf{Tnc}. We say that this estimate is optimal if there exist $c>1$ and an operator $T$ satisfying \rf{Tnc} such that
$$
\|\f(T)\|\ge\const\|\f\|_\flat.
$$
It is easy to see that in this case the completion $\cp^\flat$ of the space of polynomials in the norm 
$\|\cdot\|_\flat$ is a Banach algebra with respect to pointwise multiplication. Moreover, this Banach algebra must be an  {\it operator algebra}, i.e., it must be isomorphic to a subalgebra of all operators on Hilbert space.

Note that it was shown in \cite{Pe82} that ${\rm VMOA}\Ast H^1$ and $\big(B^0_{\be,1}\big)_+$ are Banach algebras with respect to pointwise multiplication. G. Bennett in his review MR0658618 of the paper \cite{Pe82}, published in Mathematical Reviews, showed that the same is true for the space ${\mathscr L}$. 

It turns out that the algebra $\big(B^0_{\be,1}\big)_+$ is not an operator algebra. Otherwise the norms 
$\|\cdot\|_{\mathscr L}$, $\|\cdot\|_{{\rm VMOA}\Ast H^1}$ and 
$\|\cdot\|_{B^0_{\be,1}}$ would be equivalent. It is easy to show (see \cite{Pe82}) that
this would mean that $(B^0_{\be,1})^*=B^0_{1,\be}={\frak M}H^1$. However, it is well known
(see references in \cite{Pe82}) that ${\frak M}H^1$ is a proper subset of the space $B^0_{1,\be}$.

We would like to mention here the papers \cite{Ha}, \cite{V}  and \cite{BGT}, in which the functional calculus for power bounded operators  on the Besov class $\big(B^0_{\be,1}\big)_+$ is generalized to semigroups of operators.

%
%

\subsection{Hankel matrices and tensor products of the spaces  $\bs{\ell^1}$ and $\bs{\ell^\be}$}
Let us now return to estimate \rf{cherezL}. On the right-hand side of inequality \rf{cherezL} one has to minimize 
$\|\{\g_{jk}\}_{j,k\ge0}\|_{\ell^1\widecheck\otimes\ell^1}$ over all matrices $\{\g_{jk}\}_{j,k\ge0}$ in
$\ell^1\widecheck\otimes\ell^1$ satisfying the condition 
$$
\hat\f(n)=\sum_{j+k=n}\g_{jk},\quad n\ge0.
$$
A natural question arises of whether the choice would be optimal if we put 
\bay
\label{vybor}
\g_{jk}=\frac{\hat\f(n)}{n+1}\quad\mbox{in the case when}\quad j+k=n+1.
\ey
This leads to the question of describing the Hankel matrices of class $\ell^1\widecheck\otimes\ell^1$, which was answered by the following theorem obtained in \cite{Pe82}.

\begin{thm}
\label{inGan}
Let $\f$ be a function analytic in $\dd$. Then $\G_\f\in\ell^1\widecheck\otimes\ell^1$ if and only if 
$\f\in\big(B_{\be,1}^1\big)_+$.
\end{thm}

Therefore if we make the choice \rf{vybor}, we obtain the estimate 
$$
\|\f(T)\|\le\const\|\f\|_{B^0_{\be,1}},
$$
which is not optimal as we have already mentioned earlier.

Let us return now to the averaging projection ${\mathcal P}$ onto the set of Hankel matrices, defined by \rf{proekt}. The following result obtained in \cite{Pe82} describes the image of the space $V^2$ under this projection.

\begin{thm}
\label{usprV2}
$$
{\mathcal P}V^2=\big\{\G_\f:~\f\in\big(B^0_{1,\be}\big)_+\big\}.
$$
\end{thm}

This implies easily the following fact (see \cite{Pe82}):

\begin{cor}
The projection ${\mathcal P}$ is not a bounded operator on $\ell^1\widecheck\otimes\ell^1$, 
$\ell^\be\widehat\otimes\ell^\be$ and $V^2$.
\end{cor}

\subsection{Problems 8 and 88 from the Scottish Book} In this subsection we are going to see how the result of the previous subsection allow one to solve Problems 8 and 88 from the famous Scottish Book (see \cite{Ma}). Both problems were posed by S. Mazur.

\medskip

{\bf Problem 8.} Let $c$ be the space of convergent sequences. Consider the bilinear form ${\mathscr B}$ on $c\times c$ defined by
$$
{\mathscr B}\big(\{x_n\}_{n\ge0},\{y_n\}_{n\ge0}\big)=\{z_n\}_{n\ge0},
$$
where
$$
z_n=\frac1{n+1}\sum_{k=0}^nx_ky_{n-k},\quad n\ge0.
$$
It is easy to see that ${\mathscr B}$ maps $c\times c$ to $c$. Problem 8 in the Scottish Book is whether {\it  ${\mathscr B}$ maps $c\times c$ onto the whole space $c$}.

\medskip

{\bf Problem 88.} Suppose that a Hankel matrix $\{\g_{j+k}\}$ belongs to the tensor product
$\ell^1\widecheck\otimes\ell^1$. Problem 88 of the Scottish Book is whether {\it this implies that
$$
\sum_{k\ge0}(1+k)|\g_k|<\be,
$$
that is, whether the sum of the matrix entries  must be finite}.

\medskip

In the papers \cite{KwaP} and \cite{EL} solutions to both problems were given.

It turned out, however, that it follows right away from the {\it earlier paper} \cite{Pe82} that both Problems 8 and 88 have negative solutions, see \cite{PeStud}, in which a detailed discussion is given.

\medskip

Let us proceed to the discussion of the solution to Problem 8. Consider the linear operator ${\mathscr A}$, defined on the set of matrices $Q=\{q_{jk}\}_{j,k\ge0}$ and taking matrices to sequences of complex numbers:
$$
{\mathscr A}Q=\{z_n\}_{n\ge0},\quad\mbox{where}\quad z_n\df\frac1{n+1}\sum_{j+k=n}q_{jk}.
$$

Theorem \ref{usprV2} of this survey obtained in \cite{Pe82} can be reformulated in the following way:

\begin{thm}
The operator ${\mathscr A}$ maps $V^2$ onto the space of Fourier coefficients of the functions in the Besov class $\big(B^0_{1,\be}\big)_+$.
\end{thm}

This theorem implies that
$$
{\mathscr A}\big(c\widehat\otimes c\big)\subset{\mathscr A}\big(\ell^\be\widehat\otimes \ell^\be\big)
\subset{\mathscr A}(V^2)\subset\big\{\{f_n\}_{n\ge0}:~f\in\big(B^0_{1,\be}\big)_+\big\},
$$
snd so
$$
{\mathscr B}(c\times c)\subset\big\{\{f_n\}_{n\ge0}:~f\in\big(B^0_{1,\be}\big)_+\big\}.
$$
However, it is well known that
$$
c\not\subset\big\{\{f_n\}_{n\ge0}:~f\in\big(B^0_{1,\be}\big)_+\big\},
$$
see \cite{PeStud} where a detailed discussion of such issues is given.

\medskip

Let us now proceed to Problem 88. First of all, it follows from Theorem \ref{inGan} in this survey that
$$
\{\g_{j+k}\}_{j,k\ge0}\in\ell^1\widecheck\otimes\ell^1\quad\mbox{if and only if}\quad
\sum_{n\ge0}\g_nz^n\in\big(B^1_{\be,1}\big)_+.
$$
A negative solution to Peoblem 88 follows from the following best possible estimate of the moduli of the Fourier coefficients of the functions in the Besov class $\big(B^1_{\be,1}\big)_+$:
$$
\sum_{n\ge0}2^n\left(\sum_{k=2^n}^{2^{n+1}-1}|\hat f(k)|^2
\right)^{1/2}<\be.
$$
Indeed, it is easy to see that this condition on the coefficients does not imply that 
$$
\sum_{n\ge0}n|\hat f(n)|<\be.
$$
We refer the reader to \cite{PeStud}, where a detailed discussion of such issues is given.

Note also that in  \cite{PS} certain quantitative estimates were obtained that are related to Problems 8 and 88
as well as certain related problems were posed. In the paper \cite{PeStud} it was shown that the quantitative estimates also follow from the results of \cite{Pe82}, while the problems posed in
\cite{PS} can easily be solved with the help of the results of \cite{Pe82}.

\

\section{\bf Besov classes in perturbation theory}
\setcounter{equation}{0}
\label{vozm}

\

\subsection{Functions of self-adjoint operators under their perturbation}

Let $f$ be a continuous complex-valued function on the real line. Let us consider the problem of when the inequality
\bay
\label{OLip}
\|f(A)-f(B)\|\le\const\|A-B\|
\ey
holds for arbitrary (bounded) self-adjoint operators $A$ and $B$ on Hilbert space. Functions $f$, satisfying this property are called {\it operator Lipschitz}. We denote by $\OL(\R)$ the space of operator Lipschitz functions on $\R$ and for $f\in\OL(\R)$, we put
$$
\|f\|_{\OL(\R)}\df\sup\frac{\|f(A)-f(B)\|}{\|A-B\|},
$$
where the supremum is taken over all pairs of self-adjoint operators $A$ and $B$ with bounded difference such that $A\ne B$.

It is well known (see \cite{APOL}) that if inequality \rf{OLip} holds for all bounded self-adjoint operators, then it also holds for unbounded self-adjoint operators.

Obviously, operator Lipschitz functions must be Lipschitz functions, i.e., the following inequality must be satisfied:
$$
|f(x)-f(y)|\le\const|x-y|,\quad x,~y\in\R.
$$
In the paper \cite{FaO} by Yu.B. Farforovskaya it was shown that the converse is false.

Moreover, as became clear later, operator Lipschitz functions are necessarily differentiable everywhere  (\cite{Wi}, see also \cite{APOL}) but not necessarily continuously defferentiable  (\cite{KS1}, see also \cite{APOL}). In particular, the function $x\mapsto|x|$ is not operator Lipschitz.

Next, in the papers \cite{PeFAP} and \cite{PeMD90} necessary conditions for operator Lipschitzness were obtained. In particular, if $f$ is an operator Lipschitz function, then $f$ belongs to the Besov class $B_1^1(\R)$ locally. This necessary condition also shows that the conditions of continuous differentiability and the boundedness of the derivative are not sufficient for operator Lipschitzness. We refer the reader to \cite{APOL} for other necessary conditions.

In the author's papers \cite{PeFAP} and \cite{PeMD90} it was shown that the Besov class $B_{\be,1}^1(\R)$ {\it consists of operator Lipschitz functions}. In other words, the following assertion holds:

\begin{thm}
\label{BesLip}
Let $f$ be a function in the Besov class $B_{\be,1}^1(\R)$. Then inequality {\em\rf{OLip}} holds for arbitrary self-adjoint operators $A$ and $B$.
\end{thm}

To prove Theorem \ref{BesLip}, double operator integrals were used in \cite{PeFAP} and \cite{PeMD90}. Double operator integrals are expressions of the form
\bay
\label{doi}
\int_\X\int_\Y\Phi(x,y)\,\di E_1(x)Q\,\di E_2(y),
\ey
where $E_1$ and $E_2$
are spectral measures on Hilbert space, $Q$ is a bounded linear operator on Hilbert space and $\Phi$ is a bounded measurable functions of two variables.

Double operator integrals appeared first in the paper \cite{DK}. Later, in the series of papers \cite{BS1}, \cite{BS2} and \cite{BS3} M.Sh. Birman and M.Z. Solomyak created a beautiful theory of double operator integrals. In particular, their approach allows one to define double operator integrals of the form \rf{doi} for arbitrary operators $Q$ of Hilbert--Schmidt class $\bS_2$; in this case the double operator integral must also belong to $\bS_2$ and the following inequality holds
$$
\left\|\int_\X\int_\Y\Phi(x,y)\,dE_1(x)Q\,dE_2(y)\right\|_{\bS_2}\le\|\Phi\|_{L^\be}\|Q\|_{\bS_2}.
$$
It turns out, however, that if $Q$ is an arbitrary bounded operator, then to define the double operator integral \rf{doi}, we have to impose additional assumtions on the function $\Phi$. In the paper \cite{PeFAP} the natural maximal class of such functions $\Phi$, for which the double operator integral \rf{doi} can be defined for arbitrary bounded operators $Q$ was named the {\it class of Schur multipliers} with respect to the spectral measures $E_1$ and $E_2$.

There are several descriptions of the class of Schur multipliers, see the paper \cite{PeFAP}. We mention here the following one:

\begin{thm}
A bounded measurable function  $\Phi$ is a Schur multiplier with respect to  $E_1$ and $E_2$ if and only if 
$\Phi$ belongs to the Haagerup tensor product
$L^\be_{E_1}\otimes_{\rm h}L^\be_{E_2}$, i.e., $\Phi$ admits a representation
$$
\Phi(x,y)=\sum_n\f_n(x)\psi_n(y),
$$
where $\f_n\in L^\be_{E_1}$, $\psi_n\in L^\be_{E_2}$ and
$$
\sum_n|\f_n|^2\in L^\be_{E_1}\quad\mbox{and}\quad\sum_n|\psi_n|^2\in L^\be_{E_2}.
$$
\end{thm}

Note here that if $\Phi$ is a Schur multiplier with respect to spectral measures $E_1$ and $E_2$, then 
\bay
\label{doiSp}
Q\in\bS_p,\quad1\le p<\be,\quad\Longrightarrow\quad
\int_\X\int_\Y\Phi(x,y)\,\di E_1(x)Q\,\di E_2(y)\in\bS_p.
\ey

The role of double operator integrals in obtaining operator Lipschitz estimates can be demonstrated by the following result by Birman and Solomyak \cite{BS1}. 

\begin{thm}
Let $f$ be a differentiable function on $\R$ and let $\dg f$ the corresponding divided difference, i.e., 
$$
(\dg f)(x,y)=\left\{\begin{array}{ll}\frac{f(x)-f(y)}{x-y},&x\ne y\\[.2cm]
f'(x),&x=y.
\end{array}\right.
$$
Suppose that $A$ and $B$ are self-adjoint operators such that the operator $B-A$ is bounded
$\dg f$ is a Schur multiplier with respect to the spectral measures $E_B$ and $E_A$ of the operators $B$ and $A$. Then $f(B)-f(A)$ is a bounded operator and
\bay
\label{fAB}
f(B)-f(A)=\int_\R\int_\R(\dg f)(x,y)\,\di E_B(x)(B-A)\,\di E_A(y).
\ey
\end{thm}

It follows that if the divided difference $\dg f$ is a Schur multiplier with respect to all pairs of Borel spectral measures on $\R$, then $f$ is an operator Lipschitz function. The converse is also true, see \cite{APOL} and references given there.

The following result was obtained in \cite{PeFAP}.

\begin{thm}
\label{BesOL}
Let $f$ a function in the Besov class $B_{\be,1}^1(\R)$. Then ${\mathcal D}f$ is a Schur multiplier with respect to arbitrary pairs of spectral Borel measures on $\R$, $f$ is an operator Lipschitz function and for arbitrary self-adjoint operators $A$ and $B$ with bounded difference identity {\em\rf{fAB}} is satisfied and the following inequality
holds
\bay
\label{OLBes}
\|f(A)-f(B)\|\le\const\|f\|_{B_{\be,1}^1}\|A-B\|.
\ey
\end{thm}

Note here that the condition $f\in B_{\be,1}^1(\R)$ is not necessary for operator Lipschitzness. We refer the reader to the survey \cite{APOL}, where the reader can find a detailed discussion of conditions for operator Lipschitzness.

\subsection{Operator Bernstein type inequality}

It is easy to deduce from Theorem \ref{BesOL} the following so-called operator Bernstein type inequality: let $f$ be a function of class $L^\be(\R)$ whose Fourier transform is supported in the interval $[-\s,\s]$. Then for arbitrary self-adjoint operators $A$ and $B$ with bounded difference, the following inequality holds
\bay
\label{tipBe}
\|f(A)-f(B)\|\le\const\s\|f\|_{L^\be}.
\ey
In the paper \cite{APUnb} (see also \cite{APOL}) an alternative, more elementary approach to the operator Lipschitzness of functions of class $B_{\be,1}^1(\R)$ was proposed. The first step of this approach is an independent proof of the operator Bernstein type inequality. Then the operator Lipschitzness of the functions in $B_{\be,1}^1(\R)$ was deduced from it. Moreover, in \cite{APUnb} 
it was shown that actually inequality \rf{tipBe} holds with constant 1.

In this subsection we give a proof of the operator Bernstein type inequality with constant
1 and deduce from it inequality \rf{OLBes}. We also give here a stronger operator Bernstein type inequality obtained in \cite{APUnb}, which is an operator version of an inequality obtained by S.N. Bernstein in \cite{Bern}.


Let $\s>0$. We say that an entire function $f$ is a function {\it of exponential type at most} $\s$ if for each positive number $\e$
there exists $c>0$ such that $|f(z)|\le c e^{(\s+\e)|z|}$ for all $z\in\C$.

Denote by $\E_\s$ the set of all entire functions of exponential type at most $\s$. 
Put $\E_\s^\be\df\E_\s\cap L^\be(\R)$.
It is well known that
$$
\E_\s^\be=\{f\in L^\be(\R):\supp\F f\subset[-\s,\s]\}.
$$

Note also that the space $\E_\s^\be$
coincides with the set of entire functions $f$ such that $f\in L^\be(\R)$ and
\bay
\label{fEs}
|f(z)|\le e^{\s|\im z|}\|f\|_{L^\infty(\R)},\quad z\in\C,
\ey
see, for example, \cite{Lev}, page 97.

Bernstein's inequality (see \cite{Bern}) asserts that
$$
\sup_{x\in \R}|f^\prime(x)|\le\s\sup_{x\in \R}|f(x)|
$$
for any function $f$ of class $L^\be(\R)$ whose Fourier transform has support in $[-\s,\s]$. It implies that
\bay
\label{corber}
|f(x)-f(y)|\le\s\|f\|_{L^\infty}|x-y|,\quad f\in\E_\s^\be,\quad x,~y\in\R.
\ey

Bernstein also proved in \cite{Bern} the following improvement of inequality \rf{corber}:
\bay
\label{bern2}
|f(x)-f(y)|\le\b(\s|x-y|)\|f\|_{L^\infty(\R)},\quad f\in\E_\s^\be\quad x,~y\in\R,
\ey
where
$$
\b(t)\df\left\{\begin{array}{ll}2\sin(t/2),&\text {if}\,\,\,\,0\le t\le\pi,\\[.1cm]
2,&\text {if}\,\,\,\,t>\pi.
\end{array}\right.
$$
Obviously, $\b(t)\le\min(t,2)$  for all $t\ge0$.

In the paper \cite{APUnb} the following version of the operator Bernstein type inequality was deduced from \rf{bern2}:

\begin{thm}
\label{opbern}
Let $f\in\mathscr \E_\s^\be$.  Then
\bay
\label{opnerBer}
\|f(A)-f(B)\|\le\b(\s\|A-B\|)\|f\|_{L^\infty}\le\s\|f\|_{L^\infty}\|A-B\|
\ey
for arbitrary (bounded) self-adjoint operators $A$ and $B$. In particular,
$\|f\|_{\OL(\R)}\le\lb\s\|f\|_{L^\be(\R)}$.
\end{thm}

\Pf First of all, let us observe that inequality \rf{bern2} implies the following vector version of it. 
Let $X$ be a complex Banach space. We denote by $\E_\s(X)$  the space of entire $X$-valued functions $f$ of exponential type at most $\s$, i.e., satisfying the condition:
for each positive $\e$, there exists  $c>0$ such that $\|f(z)\|_X\le c e^{(\s+\e)|z|}$ for all $z\in\C$.
Then the following inequality holds 
\bay
\label{bn}
\|f(x)-f(y)\|_X\le\b(\s|x-y|)\|f\|_{L^\infty(\R,X)}\le \s\|f\|_{L^\be(\R,X)}|x-y|,\quad x,~y\in\R.
\ey
It is easy to see that inequality \rf{bn} reduces to inequality \rf{bern2} with the help of the Hahn--Banach theorem.

Let now $A$ and $B$ be self-adjoint operators in Hilbert space $\h$. Let us show that
$$
\|f(A)-f(B)\|\le\b(\s\|A-B\|)\|f\|_{L^\infty}.
$$
Put $F(z)=f(A+z(B-A))$. Clearly, $F$ is an entire function that takes values in the space  
$\mB(\h)$ of bounded linear operators on $\h$ and $\|F(t)\|\le\|f\|_{L^\be(\R)}$ for $t\in\R$. It follows from von Neumann's inequality (see \cite{SNF}, Ch. I, \S\:8) that 
$F\in\E_{\s\|B-A\|}(\mB(\h))$. To complete the proof, it remains to apply the vector version of Bernstein's
inequality \rf{bn}  to the vector-valued function  $F$ for $x=0$ and  $y=1$.  $\bl$

\medskip

{\bf Proof of inequality \rf{OLBes} based on the operator Bernstein inequality.}
Let $f\in B_{\be,1}^1(\R)$. Consider the functions $f_n$ defined by \rf{fn}. 
Clearly, $\supp\F f_n\subset\big[-2^{n+1},2^{n+1}]$. Then by \rf{opnerBer}, we have
$$
\|f(A)-f(B)\|\le\sum_{n\in\Z}\|f_n(A)-f_n(B)\|\le\sum_{n\in\Z}2^{n+1}\|f_n\|_{L^\be}\le\const\|f\|_{B_{\be,1}^1}.
\quad\bl
$$

\medskip

{\bf Commentary.} It can be shown that if $f$ is an operator Lipschitz function on $\R$, and $A$ 
and $B$ are self-adjoint operators such that $A-B\in\bS_p$, $1\le p<\be$, then $f(A)-f(B)\in\bS_p$. However, as was shown in \cite{PoSu}, 
for $1<p<\be$ the conclusion $f(A)-f(B)\in\bS_p$ holds under a weaker assumption on $f$; it suffices to assume that $f$ a Lipschitz condition. This solved a famous problem.

\medskip

By analogy with operator Lipschitz functions, we could introduce the notion of {\it operator H\"older functions}. Let $0<\a<1$. We can say that a function $f$ on the real line $\R$ is called operator H\"older of order $\a$ if the following inequality holds
\bay
\label{Hold}
\|f(A)-f(B)\|\le\const\|A-B\|^\a
\ey
for arbitrary self-adjoint (no matter, bounded or unbounded) operators $A$ and $B$. However, it turned out that  in contrast with the case of operator Lipschitz functions, a function $f$ is operator H\"older of order $\a$ if and only if it is a H\"older function of order $\a$, i,e., belongs to $\L_\a(\R)$. This was established in \cite{APLa}. Note that in the same paper \cite{APLa} other estimates for function of operators under perturbation are found.  Another approach to the proof of inequality \rf{Hold} for functions $f$ in $\L_\a(\R)$ was found in \cite{NF}.

Note also here that in the paper \cite{APSp} the authors obtained estimates of operator differences 
$f(A)-f(B)$ in the case when $A$ and $B$ are self-adjoint operators with difference in $\bS_p$ and $f$ is a function in the H\"older class $\L_\a(\R)$. In particular, it was shown in \cite{APSp} that if $0<\a<1$, $p>1$, $f\in\L_\a(\R)$, and $A$ and $B$ are self-adjoint operators such that $A-B\in\bS_p$, then
$$
f(A)-f(B)\in\bS_{p/\a}\quad\mbox{and}\quad
\|f(A)-f(B)\|_{\bS_{p/\a}}\le\const\|f\|_{\L_\a}\|A-B\|_{\bS_p}^\a.
$$

Let us also mention here the paper \cite{McS}, in which the authors obtained estimates of the quas-norms
\lb$\|f(A)-f(B)\|_{\bS_p}$ for self-adjoint operators $A$ and $B$ in the case $0<p<1$.

\subsection{Functions of unitary operators under perturbation} The notion of operator Lipschitz functions can easily be generalized to functions on closed subsets of the complex plane $\C$. In particular, a continuous function $f$ on the unit curcle $\T$ is called {\it operator Lipschitz} if
$$
\|f(U)-f(V)\|\le\const\|U-V\|
$$
for arbitrary unitary operators $U$ and $V$ on Hilbert space.

In the paper \cite{PeFAP} the following result was obtained:

\begin{thm}
\label{unitOL}
Let $f$ be a function in the Besov class $B_{\be,1}^1(\T)$. Then $f$ is an operator Lipschitz function on $\T$ and the following inequality holds
\bay
\label{UnOL}
\|f(U)-f(V)\|\le\const\|f\|_{B_{\be,1}^1}\|U-V\|
\ey
for arbitrary unitary operators $U$ and $V$ on Hilbert space.
\end{thm}

In \cite{PeFAP} Theorem \ref{unitOL} was proved with the help of a representation of the difference of functions of the operators in terms of the double operator integral:
$$
f(U)-f(V)=\int_\T\int_\T\frac{f(\z)-f(\t)}{\z-\t}\,\di E_U(\z)(U-V)\,\di E_V(\t).
$$

Here $E_U$ and $E_V$ are the spectral measures of $U$ and $V$. More precisely, it was proved in \cite{PeFAP} that the divided difference in the integrand is a Schur multiplier and its norm in the space of Schur multipliers can be estimated in terms of the norm of $f$ in the Besov class $B_{\be,1}^1(\T)$.

As in the case of self-adjoint operators, in the paper \cite{APUnb} an alternative approach to obtain operator Lipschitz estimate \rf{UnOL} was given. That approach is also based on an operator Bernstein type inequality for functions of unitary operators:
$$
\|q(U)-q(V)\|\le(\deg q)\|U-V\|
$$
for any trigonometric polynomial $q$. 

\medskip

{\bf Commentary.} Note briefly that the problem of obtaining Lipschitz type estimates for functions of contractions was considered first in the author's paper \cite{PeSzh} (an operator is called a {\it contraction} if its norm is at most 1). In that paper with the help of double operator integrals with respect to semispectral measure it was shown that if $f$ is a function in the Besov class 
$(B_{\be,1}^1)_+$ of functions analytic in the disc $\dd$ and $T$ and $R$ are contractions on Hilbert space, then
$$
\|f(T)-f(R)\|\le\const\|f\|_{B_{\be,1}^1}\|T-R\|.
$$
We refer the reader to the paper \cite{KS4}, in which the authors consider various issues related to the behaviour of functions of contractions under perturbation.

Note also that an analogue of this result for functions of dissipative operators was obtained in  \cite{APd}.

\subsection{Functions of normal operators under their perturbation. Functions of collections of commuting operators} The notion of operator Lipschitz functions can be generalized easily to functions defined on the complex plane $\C$. Let $f$ be a continuous complex-valued function on $\C$.
We say that $f$ is an {\it operator Lipschitz function} if the inequalty
\bay
\label{OLnor}
\|f(N_1)-f(N_2)\|\le\const\|N_1-N_2\|
\ey
holds for arbitrary normal operators $N_1$ and $N_2$ with bounded difference.

As in the case of self-adjoint operators, if inequality \rf{OLnor} holds for all bounded normal operators, it also holds for unbounded normal operators, see \cite{APOL}.

It was shown in \cite{APPS} that the functions in the Besov class $B_{\be,1}^1(\R^2)$ are operator Lipschitz functions on $\C$ (here we identify naturally $\C$ with $\R^2$). 

Let us introduce the following notation. Put
$$
A_j\df\re N_j=\frac12(A_j+A_j^*),\quad B_j\df\im N_j=\frac1{2\ri},(A_j-A_j^*)
$$
$$
\quad E_j\quad\mbox{is the spectral measure of}\quad N_j,\quad j=1,\,2.
$$
In other words, $N_j=A_j+{\rm i}B_j$, $j=1,\,2$, where $A_j$ and $B_j$ are commuting self-adjoint operators. 

If $f$is a function on $\R^2$ that has partial derivatives in each variable, consider the divided differences in each variable 
$$
\big(\dg_xf\big)(z_1,z_2)\df\frac{f(x_1,y_2)-f(x_2,y_2)}{x_1-x_2},
\quad z_1,\,z_2\in\C,
$$
and
$$
\big(\dg_yf\big)(z_1,z_2)\df\frac{f(x_1,y_1)-f(x_1,y_2)}{y_1-y_2},
\quad z_1,\,z_2\in\C,
$$
where
$$
x_j\df\re z_j,\quad y_j\df\im z_j,\quad  j=1,\,2.
$$
Note that on the sets $\{(z_1,z_2):~x_1=x_2\}$ and
$\{(z_1,z_2):~y_1=y_2\}$
the divided differences $\dg_xf$ and $\dg_yf$ will be understood as the corresponding partial derivatives of $f$.

The proof of the operator Lipschitzness of function in the Besov class $B_{\be,1}^1(\R^2)$ is based on the following fact. To state it, we introduce the following notation:
$$
\E_\s^\be(\R^2)\df\big\{f\in L^\be(\R^2):~\supp\F f\subset\{\z\in\C:~|\z|\le\s\big\},
$$
where $\s$ is a positive number.

\begin{lem}
Let $f$ be a function of class $\E_\s^\be(\R^2)$.
Then $\dg_xf$ and
$\dg_yf$ are Schur multipliers for all Borel spectral measures and their norms in the space of Schur multipliers are less than or equal to
$\const\s\|f\|_{L^\be}$.
\end{lem}

\begin{thm}
\label{doin}
Let $f$ be a function in the Besov class $B_{\be,1}^1(\R^2)$. Suppose that $N_1$ and $N_2$ are normal operators such that the operator $N_1-N_2$ is bounded. Then
\begin{align*}
f(N_1)-f(N_2)&=\iint\limits_{\C^2}\big(\dg_yf\big)(z_1,z_2)\,
\di E_1(z_1)(B_1-B_2)\,\di E_2(z_2)\nonumber\\[.2cm]
&+\iint\limits_{\C^2}\big(\dg_xf\big)(z_1,z_2)\,
\di E_1(z_1)(A_1-A_2)\,\di E_2(z_2),
\end{align*}
and the following inequality holds
$$
\|f(N_1)-f(N_2)\|\le\const\|f\|_{B_{\be,1}^1}\|N_1-N_2\|,
$$
that is, $f$ is an operator Lipschitz function on $\C$.
\end{thm}

Note that the study of the behaviour of functions of normal operators under perturbation is equivalent to the problem of the behaviour functions of pairs of commuting self-adjoint operators. Naturally, we arrive at the problem of studying the behaviour of functions of $d$ commuting self-adjoint operators. Such a problem was studied in the paper \cite{NP}, in which the following result was obtained.

\begin{thm}
Let $f$ be a function in the Besov class $B_{\be,1}^1(\R^d)$. Then the following assertions hold:

{\em(i)} there exists a positive number $K$ such that
$$
\|f(A_1,\cdots,A_d)-f(B_1,\cdots,B_d)\|\le K\|f\|_{B_{\be,1}^1}
\max_{1\le j\le d}\{\|A_j-B_j\|\}
$$
for arbitrary sets $\{A_1,\cdots,A_d\}$ and $\{B_1,\cdots,B_d\}$ of commuting self-adjoint operators such that the operators $A_j-B_j$, $1\le j\le d$, are bounded;

{\em(ii)} if $\{A_1,\cdots,A_d\}$ and $\{B_1,\cdots,B_d\}$ are sets of commuting self-adjoint operators such that the operators $A_j-B_j$, $1\le j\le d$, belong to $\bS_1$, then
$f(A_1,\cdots,A_d)-f(B_1,\cdots,B_d)$ also belongs to $\bS_1$.
\end{thm}

\medskip

{\bf Commentary.} Let us mention here the papers \cite{PeCS} and \cite{APCD}, in which analogs of these results for functions of pairs of commuting contractions and functions of pairs of commuting dissipative operators were obtained.

\subsection{Operator differentiability} \label{diff}
We proceed to the question of the differentiability of the map  
\bay
\label{difvA}
K\mapsto f(A+K)-f(A).
\ey
Here $f$ is a continuous function on $\R$, $A$ is a self-adjoint operator and $K$ is a bounded self-adjoint operator. The function $f$ is called {\it operator differentiable} if this map defined on the real Banach space of bounded self-adjoint operators is differentiable for an arbitrary self-adjoint operator $A$. We refer the reader to  \cite{APOL} for a detailed discussion of operator differentiability and for an extensive bibliography; see also the papers \cite{KS2} and \cite{KS3}.

In the definition of operator differentiability given above nothing is said whether we mean Gateaux differentiability or Fr\'echet differentiability. It turns out that in this case these notions are equivalent, see \cite{APOL}.

Let us give here a sufficient condition for operator differentiability.

\begin{thm}
Let $f$ be a function in the Besov class $B_{\be,1}^1(\R)$. Then the map {\em\rf{difvA}}
is Gateaux differentiable for an arbitrary self-adjoint operator $A$. Moreover, its differential at $\0$ is the bounded linear transformer
$$
K\mapsto\int_\R\int_\R\frac{f(x)-f(y)}{x-y}\,\di E_A(x)K\,\di E_A(y).
$$
\end{thm}

We can also mention here the paper \cite{PeDifco}, in which differentiability properties of functions of contractions are studied.

\medskip

Let us proceed to operator differentiability of higher order. We start with the case when $A$ is a bounded self-adjoint operator. Then for a bounded self-adjoint operator $K$ and for a continuous function $f$ on $\R$, we   consider the problem of the existence of the $n$th order derivative of the function
\bay
\label{tf(A+tK)}
t\mapsto f(A+tK).
\ey
To obtain a formula for the $n$th order derivative of this operator-valued function, we need the notions of higher order divided differences and multiple operator integrals.

\medskip

{\bf Definition.} Let $f$ be a function on $\R$ that has $k$ derivatives. We define the {\it divided difference 
$\dg^k f$ of order} $k$ by induction in the following:
$$
\dg^0f\df f;
$$
if $k\ge1$, then
$$
(\dg^{k}f)(t_1,\cdots,t_{k+1})\df
\left\{\begin{array}{ll}
\frac{(\dg^{k-1}f)(t_1,\cdots,t_{k-1},t_k)-
(\dg^{k-1}f)(t_1,\cdots,t_{k-1},t_{k+1})}{t_{k}-t_{k+1}},&t_k\ne t_{k+1},\\[.2cm]
\frac{\partial}{\partial s}
\Big(\big(\dg^{k-1}f\big)(t_1,\cdots,t_{k-1},s)\Big)\Big|_{s=t_k},&t_k=t_{k+1}.
\end{array}\right.
$$
Note also that the result does not depend on the order of variables. We are also going to use the notation
$$
\dg f=\dg^1f.
$$

To write an expression for the higher order derivatives of operator functions, we need multiple operator integrals. There are several different approaches to the definition of multiple operator integrals. Here we are going to use an approach given in the author's paper \cite{PeOD}. Note that for other problems it becomes useful to use other definitions of multiple operator integrals, see, in particular, the definitions given in Subsection \ref{nekom} of this survey.

For simplicity, we consider the case of triple operator integrals.

Let $(\X,E)$, $(\Y,F)$ and $(\cZ,G)$ be spaces with spectral measures
$E$, $F$ and $G$ on a Hilbert space $\h$. Suppose that a function
$\psi$ on $\X\times\Y\times\Z$ belongs to the {\it integral projective tensor product} 
$L^\be(E)\hat\otimes_{\rm i}L^\be(F)\hat\otimes_{\rm i}L^\be(G)$, i.e., $\psi$ admits a representation
\bay
\label{ttp}
\psi(\l,\mu,\nu)=\int_Q f(\l,x)g(\mu,x)h(\nu,x)\,\di \s(x),
\ey
where $(Q,\s)$ is a space with a $\s$-finite measure, $f$ is a measurable function on $\X\times Q$,
$g$ is a measurable function on $\Y\times Q$, $h$ is a measurable function on $\cZ\times Q$; moreover,
\bay
\label{ner}
\int_Q\|f(\cdot,x)\|_{L^\be(E)}\|g(\cdot,x)\|_{L^\be(F)}\|h(\cdot,x)\|_{L^\be(G)}\,d\s(x)<\be.
\ey
The norm $\|\psi\|_{L^\be\hat\otimes_{\rm i}L^\be\hat\otimes_{\rm i}L^\be}$ of $\psi$ in
$L^\be(E)\hat\otimes_{\rm i}L^\be(F)\hat\otimes_{\rm i}L^\be(G)$ is defined as the infimum of the left-hand side of \rf{ner} over all representations of the form \rf{ttp}. 

Suppose now that  $T_1$ and $T_2$ are bounded linear operators on $\h$. For a function $\psi$ in
$L^\be(E)\hat\otimes_{\rm i}L^\be(F)\hat\otimes_{\rm i}L^\be(G)$ of the form \rf{ttp}, we put
\begin{align}
\label{opr}
&\int\limits_\X\int\limits_\Y\int\limits_\cZ\psi(\l,\mu,\nu)
\,\di  E(\l)T_1\,\di F(\mu)T_2\,\di G(\nu)\nonumber\\[.2cm]
\df&\int\limits_Q\left(\,\int\limits_\X f(\l,x)\,\di E(\l)\right)T_1
\left(\,\int\limits_\Y g(\mu,x)\,\di F(\mu)\right)T_2
\left(\,\int\limits_\cZ h(\nu,x)\,\di G(\nu)\right)\,\di \s(x).
\end{align}

It is well known (see, for example, \cite{PeBMS} or \cite{ACDS}) that the right-hand side of \rf{opr} does not depend on the choice of a representation \rf{ttp}; in other words, triple operator integrals are well defined.

The definition of multiple operator integrals in the general case can be done in the same way. 

The following result was obtained in the author's paper \cite{PeOD}:

\begin{thm}
\label{mraz}
Let $m$ be a natural number. Suppose that
$A$ is a bounded self-adjoint operator and $f$ is a function in the Besov class $B_{\be,1}^m(\R)$.
Then the function $\dg^m f$ belongs to the integral projective tensor product
$\underbrace{L^\be(E_A)\hat\otimes_{\rm i}\cdots\hat\otimes_{\rm i}L^\be(E_A)}_{m+1}$, 
the operator function {\em{\rf{tf(A+tK)}}} is $m$ differentiable and the following equality holds
$$
\frac{\di^m}{\di s^m}f(A+sK)\Big|_{s=0}=
m!\underbrace{\int\cdots\int}_{m+1}(\dg^{m}f)(t_1,\cdots,t_{m+1})
\,\di E_A(t_1)K\cdots K\,\di E_A(t_{m+1}).
$$
\end{thm}

In the case when $A$ is not bounded, a problem arises caused by the fact that the condition $f\in B_{\be,1}^m(\R)$ for $m>1$ does not imply that $f$ is Lipschitz on the whole line. Therefore to obtain the conclusion of the theorem in the case of unbounded operators we have to add the condition 
$f\in B_{\be,1}^1(\R)$.

\subsection{Trace class perturbations and the Lifshitz--Krein formula} In his paper \cite{Li} I.M. Lifshits while considering problems of quantum statistics and crystals theory, arrived in the following trace formula. Let $A$ and $B$ be self-adjoint operators in Hilbert space with trace class difference. Then there exists an integrable real-valued function $\xi=\xi_{\{A,B\}}$ on $\R$, for which the following trace formula holds
\bay
\label{LiKr}
\trace\big(f(B)-f(A)\big)=\int_\R f'(t)\xi(t)\,\di t
\ey
for sufficiently nice functions $f$ on $\R$. Later M.G. Krein in \cite{Kr1} gave a rigorous justification of trace formula \rf{LiKr} and showed that it holds in the case when the derivative of $f$ is the Fourier transform 
of a complex Borel measure on $\R$. M.G. Krein also posed the problem to describe the maximal class of functions $f$, for which trace formula \rf{LiKr} holds for arbitrary self-adjoint operators $A$ and $B$ with trace class difference.

In the author's paper \cite{PeMD90} the following result was obtained:

\begin{thm}
Let $A$ and $B$ be self-adjoint operators with trace difference and let $\xi$ be the spectral shift function for this pair. Then trace formula  {\em\rf{LiKr}} holds for an arbitrary function $f$ in the Besov class $B_{\be,1}^1(\R)$.
\end{thm}

However, the Besov class $B_{\be,1}^1(\R)$ is not the maximal class for which trace formula \rf{LiKr} is valid. The Krein problem was solved recently in the author's paper \cite{PeSpSd}. The solution is given by the following theorem.

\begin{thm}
Let $f$ be a continuous function on $\R$. The following are equivalent:

{\em(i)} trace formula {\em\rf{LiKr}} holds for an arbitrary pair of self-adjoint operators $A$ and $B$ with trace class difference and $\xi=\xi_{\{A,B\}}$;

{\em(ii)} $f$ is an operator Lipschitz function.
\end{thm}

{\bf Commentary.} M.G. Krein in the paper \cite{Kr2} defined an integrable real-valued spectral shift function 
$\xi$ on the circle $\T$ for pairs of unitary operators $(U,V)$ with trace class difference and obtained the trace formula  
$$
\trace\big(f(U)-f(V)\big)=\int_\T f'(\z)\xi(\z)\,\di \z
$$
for functions $f$ whose derivative has absolutely convergent Fourier series. In the paper \cite{PeFAP} this formula was extended to the case of functions $f$ in the Besov class $B_{\be,1}^1(\T)$. Later, in the paper \cite{APu} the authors described the maximal class of functions $f$, for which the trace formula holds for arbitrary unitary operators with trace class difference: as in the case of self-adjoint operators, this class coincides with the class of operator Lipschitz functions (on the circle $\T$ in this case). 

Let us briefly dwell on trace formulae for functions of contractions and functions of dissipative operators. 
The problem of finding an integrable spectral shift function for a pair of contractions with trace class difference remained open for a long time, see \cite{MNP2}, where the history of this problem is discussed in detail. The problem was completely solved in \cite{MNP1} and \cite{MNP2} by different methods; see also an earlier paper \cite{MN}, in which an additional restriction was imposed. Well then, in \cite{MNP1} and \cite{MNP2} it was proved that for an arbitrary pair of contractions $(T,R)$ on Hilbert space with trace class difference $T-R$ there exists an integrable function $\xi$ on the circle $\T$ such that 
$$
\trace\big(f(T)-f(R)\big)=\int_\T f'(\z)\xi(\z)\,\di \z
$$
for an arbitrary operator Lipschitz function $f$ analytic in the unit disc. 
In the same papers \cite{MNP1} and \cite{MNP2} trace formulae were also obtained for functions of dissipative operators.

Note that earlier partial results in this direction were obtained by H. Langer \cite{La}, V.M. Adamian--H. Neidhardt \cite{AN}, A.V. Rybkin \cite{Ryb1}, \cite{Ryb2}, \cite{Ryb3}, \cite{Ryb4} and M.G. Krein \cite{Kr87}. A more detailed discussion of the history of this problem is given in the paper \cite{MNP2}.

\subsection{Trace formulae in the case of perturbations by operators of class $\bs{\bS_m}$, 
$\bs{m\ge2}$} In L.S. Koplienko's paper \cite{Ko} the author introduced a generalized spectral shift function in the case of a perturbation of a self-adjoint operator by a self-adjoint operator of class $\bS_2$ and obtained an appropriate trace formula.

Let $A$ and $B$ be self-adjoint operators such that $K\df B-A$ is an operator in the Hilbert--Schmidt class 
$\bS_2$ and  $f$ is a sufficiently nice function on $\R$. Clearly, there is no hope to expect that the operator $f(B)-f(A)$ has to be of trace class.
The idea of Koplienko was to consider the operator
$$
f(B)-f(A)-\frac{\rm d}{{\rm d}s}\Big(f(A+sK)\Big)\Big|_{s=0},
$$
which under suitable assumptions has to be in trace class and to find a trace formula for such operators.

It was shown in \cite{Ko} that there exists a unique function $\eta$ of class $L^1(\R)$ such that
\bay
\label{ktf}
\trace\left(f(B)-f(A)-\frac{\rm d}{{\rm d}s}\Big(\f(A+sK)\Big)\Big|_{s=0}\right)=\int_\R f''(x)\eta(x)\,\di x
\ey
for arbitrary rational functions $\f$ with poles off $\R$. The function $\eta$ is called the
{\it generalized spectral shift function that corresponds to the pair} $(A,B)$. 

In the author's paper \cite{PeKN} it was shown that trace formula \rf{ktf} holds under a considerably
milder assumption on $f$:

\begin{thm}
\label{fsPSS}
Let $A$ and $B$ be self-adjoint operators such that $K\df B-A\in\bS_2$. Then trace formula
{\em\rf{ktf}} holds for an arbitrary function $f$ in the Besov class $B_{\be,1}^2(\R^2)$.
\end{thm}

Note here that the results of Koplienko were extended in \cite{Nei} to the case of functions of unitary operators;
in the author's paper \cite{PeKN} those results were improved and an analogue of the Koplienko tarce formula for unitary operators was obtained in the case when the function on the unit circle belongs to the Besov class  $B_{\be,1}^2(\T)$. 

Note also that the results of \cite{Ko} were generalized to the case of functions of contractions in \cite{PoSu2}.

In the paper \cite{PSS} the Koplienko trace formula was extended to the case of perturbations by operators of class $\bS_m$, where $m$ is an integer greater than 2. 

Let $A$ be a self-adjoint operator and let $K$ be a self-adjoint operator of class $\bS_m$. To state the main result of \cite{PSS}, we define the Taylor operator polynomial  ${\mathscr T}^{(m)}_{A,K}f$ for sufficiently smooth functions $f$ by
$$
{\mathscr T}^{(m)}_{A,K}f\df f(A+K)-f(A)
-\frac{\di}{\di t}f(A+tK)\Big|_{t=0}-
\cdots-\frac1{(m-1)!}\frac{\di^{m-1}}{\di t^{m-1}}f(A+tK)\Big|_{t=0}.
$$

In \cite{PSS} the authors established that for such operators $A$ and $K$ there exists an integrable function 
$\eta_m$ on $\R$ (a {\it spectral shift function of order} $m$ for the pair $A$ and $K$)
such that for an arbitrary function $f$ on $\R$ satisfying the assumotion
\bay
\label{prF}
\F f^{(j)}\in L^1(\R),\quad 0\le j\le m,
\ey
the following trace formula holds
\bay
\label{fslporm}
\trace\left({\mathscr T}^{(m)}_{A,K}f\right)=\int_\R f^{(m)}(x)\eta_m(x)\,\di x.
\ey

Next, it was shown in \cite{APJST} that trace fomula \rf{fslporm} also holds under a weaker condition than  \rf{prF}. Namely, the following result holds:

\begin{thm}
Let $m\ge3$. Suppose that $A$ is a self-adjoint operator and $K$ is a self-adjoint operator of class $\bS_m$. Then trace formula {\em\rf{fslporm}} holds for an arbitrary function $f$
in the Besov class $B_{\be,1}^m(\R)$.
\end{thm}

In the paper \cite{APJST} more general trace formulae were given in the case of perturbation by operators of class $\bS_m$. 
Among the results of \cite{APJST}, we mention the following one:

\begin{thm}
\label{obfsl}
Let $m$ be a natural number and let $A$ and $K$ be self-adjoint operators such that $K\in\bS_m$. Suppose that $\mu$ is a complex Borel measure on $\R$. Then there exists a complex Borel measure 
$\mu_{\{A,K\}}$ such that the following equality holds
$$
\int_\R \frac{\di^{m}}{\di t^{m}}f(A+tK)\,\di\mu(t)=(-1)^m\int_\R f^{(m)}\,\di\mu_{\{A,K\}}
$$
for an arbitrary function $f$ in the Besov class $B^m_{\be1}(\R)$. Moreover, if the measure $\mu$ is absolutely continuous (with respect to Lebesgue measure), then the measure $\mu_{\{A,K\}}$ is also absolutely continuous.
\end{thm}

The following result obtained in \cite{APJST} shows that Theorem \ref{fsPSS} is a special case of Theorem \ref{obfsl}.

\begin{thm}
Let $m$ be a positive integer and let
$A$ and $K$ be self-adjoint operators such that $K\in\bS_m$.
Consider the absolutely continuous measure $\mu$ on $\R$ defined by
$$
\di\mu(t)\df\frac{(-1)^m}{(m-1)!}(1-t)^{m-1}\chi_{[0,1]}(t)\,\di t.
$$
Then the following formula holds
$$
\trace\left({\mathscr T}^{(m)}_{A,K}f\right)=\int_\R f^{(m)}(x)\,\di\mu_{\{A,K\}}(x)
$$
for an arbitrary function $f$ in the Besov class $B_{\be1}^m(\R)$.
\end{thm}

\medskip

{\bf Commentary.} Note that the results of \cite{PSS} were extended to the case of functions of contractions in \cite{PSS2} and to the case of functions of unitary operators in \cite{PSS3}.

\subsection{Functions of pairs of noncommuting operators} \label{nekom}
For not necessarily commuting self-adjoint operators $A$ and $B$ in Hilbert space and for a function $f$ on
$\R^2$, we define the function  $f(A,B)$ of the operators $A$ and $B$ as the double operator integral 
$$
f(A,B)=\int_\R\int_\R f(s,t)\,\di E_A(s)\,\di E_B(t)=\int_\R\int_\R f(s,t)\,\di E_A(s)I\,\di E_B(t)
$$
in the case when the function $f$ is a Schur multiplier with respect to the spectral measures
$E_A$ and $E_B$ of the operators $A$ and $B$.

In the case when the operators $A$ and $B$ are bounded, the functions $f$ in the Besov class 
$B_{\be,1}^1(\R)$ are Schur multipliers with respect to the spectral measures 
$E_A$ and $E_B$ (see \cite{ANP}), and so for such functions $f$ we can define the functions $f(A,B)$ of the operators $A$ and $B$.

We proceed now to the problem of the behaviour of functions of pairs of noncommuting self-adjoint operators under perturbation. In other words, we would like to estimate the norms of the difference 
$f(A_2,B_2)-f(A_1,B_1)$ in terms of the norms of $A_2-A_1$ and $B_2-B_1$. For this purpose we need triple operator integrals. Moreover, we need a definition of triple operator integrals, which differs from the one given in the subsection \ref{diff}.

Recall that triple operator integrals are expressions of the form
\bay
\label{tropin}
\int_\X\int_\Y\int_\Z\Psi(x,y,z)\,\di E_1(x)T\,dE_2(y)R\,\di E_3(z),
\ey
where $E_1$, $E_2$ and $E_3$ are spectral measure on Hilbert space, and $T$ and $R$ are bounded linear operators.

In the paper \cite{JTT} triple operator integrals were defined for functions $\Psi$ that belong to the {\it Haagerup tensor product} 
$L^\be_{E_1}\otimes_{\rm h}L^\be_{E_2}\otimes_{\rm h}L^\be_{E_3}$, which, by definition, consists of function
 $\Psi$ that admit a representation
\bay
\label{predPsi}
\Psi(x,y,z)=\sum_{j,k\ge0}\a_j(x)\b_{jk}(y)\g_k(z),
\ey
where the functions $\a_j$, $\b_{jk}$ and $\g_k$ satisfy the requirements
$$
\sum_{j\ge0}|\a_j|^2\in L^\be_{E_1},\qquad
\quad\sum_{k\ge0}|\g_k|^2\in L^\be_{E_3}
$$
and the function
$$
y\mapsto\big\|\{\b_{jk}(y)\}_{j,k}\big\|_{\mathcal B}
$$
belongs to $L^\be_{E_2}$. Here the symbol ${\mathcal B}$ denotes the class of infinite matrices that determine bounded bounded linear operators on $\ell^2$. Herewith the norm of a matrix in 
${\mathcal B}$ is defined as the norm of the corresponding operator from $\ell^2$ to $\ell^2$.

The norm of $\Psi$ in the Haagerup tensor product 
$L^\be_{E_1}\otimes_{\rm h}L^\be_{E_2}\otimes_{\rm h}L^\be_{E_3}$ is defined as the infimum of the products 
$$
\|\{\a_j\}_{j\ge0}\|_{L^\be_{E_1}(\ell^2)}\|\{\b_{jk}\}_{j,k\ge0}\|_{L^\be_{E_2}({\mathcal B})}
\|\{\g_k\}_{k\ge0}\|_{L^\be_{E_3}(\ell^2)}
$$
over all representations of $\Psi$ of the form \rf{predPsi}.

For such functions $\Psi$, the triple operator integral \rf{tropin} is defined as the product of operator matrices 
$$
\left(\begin{matrix}A_0&A_1&A_2&\cdots&\end{matrix}\right)
\left(\begin{matrix}B_{00}&B_{01}&B_{02}&\cdots&\\[.2cm]
B_{10}&B_{11}&B_{12}&\cdots\\[.2cm]
B_{20}&B_{21}&B_{22}&\cdots\\[.2cm]
\vdots&\vdots&\vdots&\ddots
\end{matrix}\right)
\left(\begin{matrix}G_0\\[.2cm]G_1\\[.2cm]G_2\\[.2cm]\vdots
\end{matrix}\right),
$$
where
$$
A_j\df\int\a_j\,\di E_1,\quad B_{jk}\df\int\b_{jk}\,\di E_2\quad\mbox{and}\quad G_k\df\int\g_k\,\di E_3.
$$

Herewith the following inequalities hold:
$$
\left\|\iiint\Psi\,\di E_1T\,\di E_2R\,\di E_3\right\|\le\|\Psi\|_{L^\be\otimes_{\rm h}L^\be\otimes_{\rm h}L^\be}
\|T\|\cdot\|R\|
$$
(see \cite{JTT}) and 
\bay
\label{<be}
\left\|\iiint\Psi\,\di E_1T\,\di E_2R\,\di E_3\right\|_{\bS_r}
\le\|\Psi\|_{L^\be\otimes_{\rm h}L^\be\otimes_{\rm h}L^\be}
\|T\|_{\bS_p}\cdot\|R\|_{\bS_q},\quad \frac1r=\frac1p+\frac1q,
\ey
in the case when $2\le p\le\be$ and $2\le q\le\be$ (see  \cite{APMOI}).

It turned out, however, that to obtain Lipschitz type estimates for functions of two noncommuting self-adjoint operators under perturbation,
the definition of triple operator integrals in the case when the integrand $\Psi$ belongs to the Haagerup tensor of $L^\be$ spaces does not help. 

In the paper \cite{ANP} the following Lipschitz type estimate for functions of pairs of noncom\-muting self-adjoint operators was obtained:

Let $1\le p\le2$. Then there exists a positive number $C$ such that
\bay
\label{Lip2nek}
\|f(A_1,B_1)-f(A_2,B_2)\|_{\bS_p}\le C\|f\|_{B_{\be,1}^1}
\max\big\{\|A_1-A_2\|_{\bS_p},\|B_1-B_2\|_{\bS_p}\big\}
\ey
for an arbitrary function $f$ in the Besov class $B_{\be,1}^1(\R^2)$ and for arbitrary bounded (not necessarily commuting) self-adjoint operators $A_1$, $A_2$, $B_1$ and $B_2$ satisfying the assumptions $A_2-A_1\in\bS_p$ and $B_2-B_1\in\bS_p$.

It is easy to deduce from the definition of Besov classes given in \S\;\ref{Besov} that to \lb prove inequality \rf{Lip2nek}, it suffices to establish that if $f\in L^\be(\R^2)$ and
\lb$\supp\F f\subset\{\xi\in\R^2:~\|\xi\|\le1\}$, then for $1\le p\le2$ the following inequality holds
\bay
\label{dlyas1}
\|f(A_1,B_1)-f(A_2,B_2)\|_{\bS_p}\le\const\|f\|_{L^\be}
\max\big\{\|A_1-A_2\|_{\bS_p},\|B_1-B_2\|_{\bS_p}\big\}
\ey
for arbitrary bounded self-adjoint operators $A_1$, $A_2$, $B_1$ and $B_2$ such that $A_2-A_1\in\bS_p$ and $B_2-B_1\in\bS_p$.

To establish inequality \rf{dlyas1}, the authors of the paper \cite{ANP} used the following representation of the operator difference $f(A_1,B_1)-f(A_2,B_2)$ in terms of triple operator integrals
\begin{align}
\label{f(A1B1)-f(A_2B_2)}
f(A_1,B_1)-&f(A_2,B_2)=
\iiint\big(\dg^{[1]}f\big)(x_1,x_2,y)
\,\di E_{A_1}(x_1)(A_1-A_2)\,\di E_{A_2}(x_2)\,\di E_{B_1}(y),\nonumber\\[.2cm]
+&\iiint\big(\dg^{[2]}f\big)(x,y_1,y_2)
\,\di E_{A_2}(x)\,\di E_{B_1}(y_1)(B_1-B_2)\,\di E_{B_2}(y_2),
\end{align}
where
\bay
\label{D1}
\big(\dg^{[1]}f\big)(x_1,x_2,y)\df\frac{f(x_1,y)-f(x_2,y)}{x_1-x_2}
\ey
and
\bay
\label{D2}
\big(\dg^{[2]}f\big)(x,y_1,y_2)\df\frac{f(x,y_1)-f(x,y_2)}{y_1-y_2}.
\ey
Certainly, we have to clarify, in which sense we understand triple operator integrals.

It turns out, (see \cite{ANP}) that under the condition $f\in\E_\1^\be(\R^2)$,
the divided differences $\big(\dg^{[1]}f\big)$ and $\big(\dg^{[2]}f\big)$ admit representations
\bay
\label{D1}
\big(\dg^{[1]}f\big)(x_1,x_2,y)=
\sum_{j,k\in\Z}\frac{\sin(x_1-j\pi)}{x_1-j\pi}\cdot\frac{\sin(x_2-k\pi)}{x_2-k\pi}
\cdot\frac{f(j\pi,y)-f(k\pi,y)}{j\pi-k\pi},
\ey
and
\bay
\label{D2}
\big(\dg^{[2]}f\big)(x,y_1,y_2)=\sum_{j,k\in\Z}
\frac{f(j\pi,x)-f(k\pi,x)}{j\pi-k\pi}\cdot
\frac{\sin(y_1-j\pi)}{y_1-j\pi}\cdot\frac{\sin(y_2-k\pi)}{y_2-k\pi},
\ey
also, in the case $j=k$ we assume that
$$
\frac{f(j\pi,y)-f(k\pi,y)}{j\pi-k\pi}
=\frac{\partial f(x,y)}{\partial x}\Big|_{(j\pi,y)}.
$$
It is well known that 
$$
\sum_{j\in\Z}\frac{\sin^2(x_1-j\pi)}{(x_1-j\pi)^2}
=\sum_{k\in\Z}\frac{\sin^2(x_2-k\pi)}{(x_2-k\pi)^2}=1,
\quad x_1~x_2\in\R,
$$
and
$$
\sup_{y\in\R}\left\|\left\{\frac{f(j\pi,y)-f(k\pi,y)}{j\pi-k\pi}
\right\}_{j,k\in\Z}\right\|_{\mB}\le\const\|f\|_{L^\be(\R)}.
$$

Similar facts can easily be stated for functions of class $\E_\s^\be(\R^2)$ for an arbitrary positive number $\s$.

It would seem the representations of the divided differences $\dg^{[1]}$ and $\dg^{[2]}$ in the form \rf{D1}
and \rf{D2} remind a representation of a function of three variables in the form of the Haagerup tensor product. However, upon closer examination, it turns out that this is not the case. The thing is that
in the definition of the Haagerup tensor product of three $L^\be$ spaces, we see that in representation \rf{predPsi} the matrix-valued function $y\mapsto\{\b_{jk}(y)\}_{j,k}$ is a function of the second variable while in representations \rf{D1} and \rf{D2} the-matrix-valued functions are functions of the first and the third variable. What is more, this distinction is very essential.

Therefore the need to define triple operator integrals for the functions $\dg^{[1]}$ and $\dg^{[2]}$ led in the paper \cite{ANP} to notions of Haagerup like tensor products (see also \cite{APMOI}).

\medskip

{\bf Definition 1.} We say that a bounded measurable function $\Psi$ on $\X\times\Y\times Z$
belongs to the {\it Haagerup-like tensor product of the first kind} 
$L^\be(E_1)\!\otimes_{\rm h}\!L^\be(E_2)\!\otimes^{\rm h}\!L^\be(E_3)$ if $\Psi$ admits a representation of the form
\bay
\label{yaH}
\Psi(x,y,z)=\sum_{j,k\ge0}\a_j(x)\b_{k}(y)\g_{jk}(z),\quad x\in\X,\quad y\in\Y,\quad z\in\Z,
\ey
where $\{\a_j\}_{j\ge0},~\{\b_k\}_{k\ge0}\in L^\be(\ell^2)$ and 
$\{\g_{jk}\}_{j,k\ge0}\in L^\be(\mB)$. Put 
$$
\|\Psi\|_{L^\be\otimes_{\rm h}\!L^\be\otimes^{\rm h}\!L^\be}
\df\inf\big\|\{\a_j\}_{j\ge0}\big\|_{L^\be(\ell^2)}
\big\|\{\b_k\}_{k\ge0}\big\|_{L^\be(\ell^2)}
\big\|\{\g_{jk}\}_{j,k\ge0}\big\|_{L^\be(\mB)},
$$
where the infimum is taken over all representations of the form \rf{yaH}.

\medskip

Now it is turn to define triple operator integrals whose integrands belong to 
$L^\be(E_1)\!\otimes_{\rm h}\!L^\be(E_2)\!\otimes^{\rm h}\!L^\be(E_3)$.

\medskip

{\bf Definition 2.}
Let $1\le p\le2$. For functions 
$\Psi$ in $L^\be(E_1)\!\otimes_{\rm h}\!L^\be(E_2)\!\otimes^{\rm h}\!L^\be(E_3)$, for a bounded linear operator $R$ and for an operator $T$ of class $\bS_p$, we define the triple operator integral
\bay
\label{WHft}
W=\iiint\Psi(x_1,x_2,x_3)\,\di E_1(x_1)T\,\di E_2(x_2)R\,\di E_3(x_3)
\ey
as the following continuous linear functional on $\bS_{p'}$,
$1/p+1/p'=1$ (on the class of compact operators in the case $p=1$):
\bay
\label{fko}
Q\mapsto
\trace\left(\left(
\iiint
\Psi(x_1,x_2,x_3)\,\di E_2(x_2)R\,\di E_3(x_3)Q\,\di E_1(x_1)
\right)T\right).
\ey

\medskip

Then the triple operator integral in \rf{fko} is well defined because the function
$$
(x_2,x_3,x_1)\mapsto\Psi(x_1,x_2,x_3)
$$ 
belongs to the Haagerup tensor product
$L^\be(E_2)\!\otimes_{\rm h}\!L^\be(E_3)\!\otimes_{\rm h}\!L^\be(E_1)$. It is easy to deduce from 
\rf{<be} that
$$
\|W\|_{\bS_p}\le\|\Psi\|_{L^\be\otimes_{\rm h}\!L^\be\otimes^{\rm h}\!L^\be}
\|T\|_{\bS_p}\|R\|,\quad1\le p\le2.
$$

Let us now define triple operator integrals in the case of integrands in the Haagerup like tensor product of the second kind. 

\medskip

{\bf Definition 3.} We say that a bounded measurable function $\Psi$ on $\X\times\Y\times Z$
belongs to the {\it Haagerup like tensor product of the second kind}  
$L^\be(E_1)\!\otimes^{\rm h}\!L^\be(E_2)\!\otimes_{\rm h}\!L^\be(E_3)$ if
$\Psi$ admits a representation
\bay
\label{preds}
\Psi(x_1,x_2,x_3)=\sum_{j,k\ge0}\a_{jk}(x_1)\b_{j}(x_2)\g_k(x_3),
\ey
where $\{\b_j\}_{j\ge0},~\{\g_k\}_{k\ge0}\in L^\be(\ell^2)$, 
$\{\a_{jk}\}_{j,k\ge0}\in L^\be(\mB)$. By the norm of $\Psi$ in 
$L^\be\otimes^{\rm h}\!L^\be\otimes_{\rm h}\!L^\be$ we understand
$$
\|\Psi\|_{L^\be\otimes^{\rm h}\!L^\be\otimes_{\rm h}\!L^\be}
\df\inf\big\|\{\a_j\}_{j\ge0}\big\|_{L^\be(\ell^2)}
\big\|\{\b_k\}_{k\ge0}\big\|_{L^\be(\ell^2)}
\big\|\{\g_{jk}\}_{j,k\ge0}\big\|_{L^\be(\mB)},
$$
where the infimum is taken over all representations of the form \rf{preds}.

\medskip

{\bf Definition 4.}
Suppose now that
$\Psi\in L^\be(E_1)\!\otimes^{\rm h}\!L^\be(E_2)\!\otimes_{\rm h}\!L^\be(E_3)$,
$T$ is a bounded linear operator and $R\in\bS_p$, $1\le p\le2$. Then the continuous linear functional 
$$
Q\mapsto
\trace\left(\left(
\iiint\Psi(x_1,x_2,x_3)\,\di E_3(x_3)Q\,\di E_1(x_1)T\,\di E_2(x_2)
\right)R\right)
$$
on the space $\bS_{p'}$ (on the space of compact operators in the case $p=1$) 
determines an operator $W$ of class $\bS_p$ and we call this operator the triple operator integral 
\bay
\label{WHst}
W=\iiint\Psi(x_1,x_2,x_3)\,\di E_1(x_1)T\,\di E_2(x_2)R\,\di E_3(x_3).
\ey

Herewith the following inequality holds:
$$
\|W\|_{\bS_p}\le
\|\Psi\|_{L^\be\otimes^{\rm h}\!L^\be\otimes_{\rm h}\!L^\be}
\|T\|\cdot\|R\|_{\bS_p}.
$$

These definitions of triple operator integrals were given in the paper \cite{ANP}; we also refer the reader to the paper \cite{APMOI}, where triple operator integrals were defined in a more general situation.

Now we are ready to state the theorem on Lipschitz type estimates for functions of noncommuting self-adjoint operators under perturbation.

\begin{thm}
\label{vozfotpar}
Let $f$ be a function in the Besov class $B_{\be,1}^1(\R^2)$ and let $1\le p\le 2$. Suppose that
$(A_1,B_1)$ are $(A_2,B_2)$ pairs of bounded not necessarily commuting self-adjoint operators such that   $A_2-A_1\in\bS_p$ and $B_2-B_1\in\bS_p$.
Then the divided differences $\big(\dg^{[1]}f\big)$ and $\big(\dg^{[2]}f\big)$ defined by 
{\em\rf{D1}} and {\em\rf{D2}} belong to the Haagerup like tensor products 
$L^\be(E_1)\!\otimes_{\rm h}\!L^\be(E_2)\!\otimes^{\rm h}\!L^\be(E_3)$ and
$L^\be(E_1)\!\otimes^{\rm h}\!L^\be(E_2)\!\otimes_{\rm h}\!L^\be(E_3)$ for arbitrary Borel spectral measures $E_1$, $E_2$ and $E_3$; moreover, representation {\em\rf{f(A1B1)-f(A_2B_2)}} holds, in which the first triple operator integral should be understood in the sense of Definition {\em2} while the second triple operator integral should be understood in the sense of Definition {\em4}.

Moreover, the following Lipschitz type estimate holds 
\bay
\label{Liptip}
\|f(A_1,B_1)-f(A_2,B_2)\|_{\bS_p}\le\const\|f\|_{B_{\be,1}^1}
\max\big\{\|A_1-B_1\|_{\bS_p},\|A_2-B_2\|_{\bS_p}\big\}.
\ey
\end{thm}

Theorem \ref{vozfotpar} was established in \cite{ANP}. In the same paper it was shown that for $p>2$ as well as in the operator norm, such a Lipschitz estimate does not hold.

Note also that it was shown in \cite{Petri} that such Lipschitz type estimates cannot be generalized to the case of functions of three noncommuting self-adjoint operators.

\subsection{Functions of pairs of unbounded noncommuting operators}
We proceed now to the case of pairs of unbounded not necessarily commuting self-adjoint operators.

Let $1\le p\le 2$ and let $f$ a function in the Besov class $B_{\be,1}^1(\R)$. Suppose now that
$(A_1,B_1)$ and $(A_2,B_2)$ are pairs of unbounded not necessarily commuting self-adjoint operators such that  $A_2-A_1\in\bS_p$ and $B_2-B_1\in\bS_p$.

Note, first of all, that both triple operator integrals on the right-hand side of representation
\rf{f(A1B1)-f(A_2B_2)} are still well defined as in the case of pairs of bounded operators and determine operators of class $\bS_p$ while their norms in
$\bS_p$ can still be estimated in terms of the right-hand side of inequality \rf{Liptip}.

The problem is how to define functions of operators $f(A_1,B_1)$ and $f(A_2,B_2)$ as unbounded densely defined operators. It is to this problem the papers 
\cite{APDok} and \cite{APne} are devoted.

In Subsection \ref{nekom} we have defined the functions $f(A,B)$ of not necessarily commuting self-adjoint operators in the case when $f\in L^\be_{E_A}\otimes_{\rm h}L^\be_{E_B}$.
Let now $f$ be a function of two variables and let $f_\sharp$ be the function defined by  
$f_\sharp(s,t)\df(1-\ri t)^{-1}f(s,t)$.
Assume that $f_\sharp\in L^\be_{E_A}\otimes_{\rm h}L^\be_{E_B}$. 
We define the operator $f(A,B)$ by
$$
f(A,B)\df f_\sharp(A,B)(I-\ri B)=\left(\,\,\,\iint\limits_{\R\times\R}f_\sharp(s,t)\,\di E_A(s)\,\di E_B(t)\right)(I-\ri B).
$$
Then $f(A,B)$ is a densely defined operator whose domain coincides with the domain ${\rm D}(B)$ of $B$. It is not necessarily bounded but the operator $f(A,B)(I-\ri B)^{-1}$ is bounded.

It was established in \cite{APdi} and \cite{APne} that if $f\in\E_\s^\be(\R^2)$, $\s>0$, then
 $f_\sharp\in L^\be_{E_1}\!\otimes_{\rm h}\! L^\be_{E_2}$ for arbitrary Borel spectral measure $E_1$ and $E_2$. Therefore, if $f\in\E_\s^\be(\R^2)$, then
the operator $f_\sharp(A,B)$ is bounded while $f(A,B)$ is a not necessarily bounded densely defined operator with domain ${\rm D}(B)$. Herewith the following inequality holds
$$
\|f_\sharp\|_{L^\be_{E_1}\otimes_{\rm h}L^\be_{E_2}}\le\const(1+\s)\|f\|_{L^\be(\R^2)}.
$$
This easily implies that if
 $f$ is a function in the inhomogeneous Besov class $\Bs(\R^2)$, then 
 $f_\sharp\in L^\be_{E_1}\!\otimes_{\rm h}\! L^\be_{E_2}$
and
$$
\|f_\sharp\|_{L^\be_{E_1}\otimes_{\rm h}L^\be_{E_2}}\le\const\|f\|_{\Bs}
$$
for arbitrary spectral Borel measures $E_1$ and $E_2$ (see \cite{APDok} and \cite{APne}).

Thus, for an arbitrary pair $(A,B)$ of not necessarily bounded self-adjoint operators and for a function $f$ in the inhomogeneous Besov class $\Bs(\R^2)$, we can define the operator $f(A,B)$ as a densely defined operator. 
Herewith the following result holds
(see  \cite{APDok} and \cite{APne}):

\begin{thm}
Let $f\in\Bs(\R^2)$ and let $p\in[1,2]$. Suppose that $(A_1,B_1)$ and $(A_2,B_2)$ are pairs of not necessarily bounded and not necessarily commuting self-adjoint operators such that
$A_2-A_1\in\bS_p$ and $B_2-B_1\in\bS_p$. Then
$$
\|f(A_1,B_1)-f(A_2,B_2)\|_{\bS_p}\le\const\|f\|_{\Bs}
\max\big\{\|A_1-A_2\|_{\bS_p},\|B_1-B_2\|_{\bS_p}\big\}.
$$
\end{thm}

Note here that an analogue of this result for functions of pairs of noncommuting dissipative operators was obtained in \cite{APdi}.

\subsection{Functions of pairs of almost commuting operators} Let $A$ and $B$ be bounded self-adjoint operators. They are said to be {\it almost commuting} if their 
 {\it commutator} $[A,B]\df AB-BA$ is a trace class operator.

A polynomial calculus for a pair $(A,B)$ of almost commuting self-adjoint operators is defined in the following way. If $\f$ is a polynomial of the form $\f(s,t)=\sum_{j,k}a_{jk}s^jt^k$, then the operator $\f(A,B)$ is defined by
$$
\f(A,B)\df\sum_{j,k}a_{jk}A^jB^k.
$$
It is easy to verify that if $\f$ are $\psi$ polynomials of two real variables, then  
$[\f(A,B),\psi(A,B)]\lb\in\bS_1$.

In the paper \cite{HH} by Helton and Howe the authors obtained a trace formula for commutators of polynomials of almost commuting self-adjoint operators. Let $A$ and $B$ be almost commuting self-adjoint operators and let $\f$ and $\psi$ be polynomials of two real variables. Then
\bay
\label{HeHo}
\trace[\f(A,B),\psi(A,B)]=
-\ri\int_{\R^2}\left(\frac{\partial\f}{\partial x}\frac{\partial\psi}{\partial y}-
\frac{\partial\f}{\partial y}\frac{\partial\psi}{\partial x}\right)\di P,
\ey
where $P$ is the real-valued measure with compact support, which is determined by the pair $(A,B)$. 
It was shown in \cite{Pin2} that the measure $P$ is absolutely continuous with respect to planar Lebesgue measure and
$$
\di P(x,y)=\frac1{2\pi}g(x,y)\,\di\m_2(x,y),
$$
where $g$ is the so-called {\it Pincus principal function} introduced in \cite{Pin1}. 

The polynomial calculus was extended in  \cite{CP} to the calculus on the class of functions  $\f=\F\o$ that are Fourier transforms of complex Borel measure 
$\o$ on $\R^2$ such that
$$
\int_{\R^2}(1+|t|)(1+|s|)\,\di |\o|(s,t)<\be;
$$
moreover the Helton--Howe trace formula \rf{HeHo} was extended to the case of such functions.

Next, in the author's paper \cite{PeAC} such results were obtained for a yet broader class of functions. Moreover, in the same paper \cite{PeAC} it was shown that it is impossible to extend such a functional calculus  to the class of all continuously differentiable functions on  $\R^2$ that would satisfy certain natural requirements  and  $[\f(A,B),\psi(A,B)]\in\bS_1$ for arbitrary continuously differentiable functions $\f$ and $\psi$.

Most advanced results in this direction were obtained in the paper \cite{APKomu}, in which the desired functional calculus is constructed on the Besov class  $B_{\be,1}^1(\R^2)$ and it is shown that this functional calculus is almost multiplicative. 

Let us proceed now to the main results of 
the paper \cite{APKomu}. Note that these results use triple operator integrals for Haagerup-like tensor products of the $L^\be$ spaces that were introduced in \cite{ANP} and the estimates of the trace norms of such operator integrals that were found in \cite{ANP}, see Subsection \ref{nekom}.

\begin{thm}
\label{komut}
Let $A$ and $B$ be self-adjoit operators and let $Q$ be a bounded linear operator such  $[A,Q]\in\bS_1$ and 
$[B,Q]\in\bS_1$. 
Suppose that $\f$ is a function in the Besov class $B_{\be,1}^1(\R^2)$.
Then $[\f(A,B),Q\big]\in\bS_1$,
\begin{align*}
\label{komQ}
\big[\f(A,B),Q\big]&=
\iiint\frac{\f(x,y_1)-\f(x,y_2)}{y_1-y_2}\,\di E_A(x)\,\di E_B(y_1)[B,Q]\,\di E_B(y_2)\nonumber
\\[.2cm]
&+
\iiint\frac{\f(x_1,y)-\f(x_2,y)}{x_1-x_2}\,\di E_A(x_1)[A,Q]\,\di E_A(x_2)\,\di E_B(y)
\end{align*}
and
$$
\label{ABQ}
\big\|[\f(A,B),Q\big]\big\|_{\bS_1}
\le\const\|\f\|_{B_{\be,1}^1(\R^2)}\big(\big\|[A,Q]\big\|_{\bS_1}+
\big\|[B,Q]\big\|_{\bS_1}\big).
$$
\end{thm}

\begin{thm}
\label{glav}
Let $A$ and $B$ be almost commuting self-adjoint operators and let $\f$ and $\psi$ be functions in the Besov $B_{\be,1}^1(\R^2)$. Then
\begin{align*}
\!\!\!\!\!\big[\f(A,B),\psi(A,B)\big]\!&=
\iiint
\frac{\f(x,y_1)-\f(x,y_2)}{y_1-y_2}\,\di E_A(x)\,\di E_B(y_1)[B,\psi(A,B)]\,\di E_B(y_2)
\\[.2cm]
&+\!\!
\iiint
\frac{\f(x_1,y)-\f(x_2,y)}{x_1-x_2}\di E_A(x_1)[A,\psi(A,B)]\di E_A(x_2)\di E_B(y)
\end{align*}
and
$$
\big\|[\f(A,B),\psi(A,B)\big]\big\|_{\bS_1}
\le\const\|\f\|_{B_{\be,1}^1(\R^2)}\|\psi\|_{B_{\be,1}^1(\R^2)}
\big\|[A,B]\big\|_{\bS_1}.
$$
\end{thm}

Now we can extend the Helton--Howe trace formula to the case of functions in the Besov class 
$B_{\be,1}^1(\R^2)$.

\begin{thm}
Let $A$ and $B$ be almost commuting self-adjoint operators and let $\f$ and $\psi$ be functions in the Besov class $B_{\be,1}^1(\R^2)$. Then
\begin{align*}
\trace\big({\rm i}\big(\f(A,B)\psi(A,B)&-\psi(A,B)\f(A,B)\big)\big)\nonumber\\[.2cm]
&=\frac{1}{2\pi}
\iint_{\R^2}\left(\frac{\partial\f}{\partial x}\frac{\partial\psi}{\partial y}-
\frac{\partial\f}{\partial y}\frac{\partial\psi}{\partial x}\right)g(x,y)\,\di x\,\di y,
\end{align*}
where $g$ is the Pincus principal function for the pair $(A,B)$.
\end{thm}

\

\renewcommand\l{\pl}

\

\
\noindent

\begin{footnotesize}
\begin{tabular}{p{8cm}}
St.Petersburg State University\\
Universitetskaya nab., 7/9 \\
199034 St.Petersburg\\
Russia\\

\\

St.Petersburg Department\\
Steklov Institute of Mathematics \\
Russian Academy of Sciences\\
Fontanka 27, 191023 St.Petersburg\\
Russia\\

\\

People's Friendship University\\
of Russia (RUDN University)\\
6 Miklukho-Maklaya St, Moscow,\\
6117198, Russian Federation

\end{tabular}

\end{footnotesize}

\end{document}

%% file: classstartscr.tex
\newcommand{\lb}{\linebreak}

\let\pl = \l

\renewcommand{\a}{\alpha}
\renewcommand{\b}{\beta}
\newcommand{\g}{\gamma}

\newcommand{\e}{\varepsilon}
\newcommand{\vk}{\varkappa}
\newcommand{\z}{\zeta}

\renewcommand{\l}{\lambda}

\newcommand{\s}{\sigma}
\renewcommand{\t}{\tau}

\newcommand{\f}{\varphi}
\renewcommand{\o}{\omega}

\newcommand{\G}{\Gamma}
\newcommand{\D}{\Delta}
\renewcommand{\L}{\Lambda}

\renewcommand{\O}{\Omega}

\newcommand{\E}{{\mathscr E}}

\newcommand{\F}{{\mathscr F}}

\newcommand{\h}{{\mathscr H}}

\newcommand{\cp}{{\mathscr P}}

\newcommand{\X}{{\mathscr X}}
\newcommand{\Y}{{\mathscr Y}}

\newcommand{\1}{{\bf 1}}

\newcommand{\C}{{\Bbb C}}
\newcommand{\T}{{\Bbb T}}
\newcommand{\pp}{{\Bbb P}}
\newcommand{\dd}{{\Bbb D}}
\newcommand{\R}{{\Bbb R}}
\newcommand{\Z}{{\Bbb Z}}

\newcommand{\0}{{\boldsymbol{0}}}

\newcommand{\bs}{\boldsymbol}

\newcommand{\m}{{\boldsymbol m}}
\newcommand{\bS}{{\boldsymbol S}}

\newcommand{\rf}[1]{(\ref{#1})}

\newcommand{\df}{\stackrel{\mathrm{def}}{=}}
\newcommand{\dist}{\operatorname{dist}}

\newcommand{\re}{\operatorname{Re}}
\newcommand{\spn}{\operatorname{span}}
\newcommand{\supp}{\operatorname{supp}}
\newcommand{\clos}{\operatorname{clos}}

\newcommand{\trace}{\operatorname{trace}}
\newcommand{\rank}{\operatorname{rank}}
\newcommand{\const}{\operatorname{const}}

\newcommand{\eeq}{\end{equation}}
\newcommand{\beq}{\begin{equation}}
\newcommand{\bay}{\begin{eqnarray}}
\newcommand{\ba}{\begin{align*}}
\newcommand{\ea}{\end{align*}}
\newcommand{\ey}{\end{eqnarray}}
\newcommand{\bey}{\begin{eqnarray*}}
\newcommand{\eey}{\end{eqnarray*}}

\newcommand{\be}{\infty}

\newcommand{\bl}{\blacksquare}

\newcommand{\Pf}{{\bf Proof. }}
\newcommand{\im}{\operatorname{Im}}
\renewcommand{\re}{\operatorname{Re}}
\newcommand{\ov}{\overline}

\newtheorem{thm}{\hspace{\parindent}Theorem}[section]

\newtheorem{cor}[thm]{\hspace{\parindent}Corollary}
\newtheorem{lem}[thm]{\hspace{\parindent}Lemma}

%% file: operator_Besov_eng.bbl
\begin{thebibliography}{99}
\label{bibl}

\medskip

\bibitem{AAK} {\sc V.M. Adamyan, D.Z. Arov,} and {\sc M.G. Krein}, {\em Analytic properties of Schmidt pairs for a Hankel operator and
the generalized Schur--Takagi problem}, {\em Mat. Sb.} {\bf 86} (1971), 34--75 (Russian);
English transl. {\em Math. USSR Sb.} {\bf15} (1971), 31--73.

\bibitem{AN} {\sc V.M. Adamjan} and {\sc H. Neidhardt}, 
{\em On the summability of the spectral shift function for pair of
 contractions and dissipative operators},
J. Operator Th. {\bf24} (1990), 187--205.

\bibitem{ANP} {\sc A.B. Aleksandrov, F.L. Nazarov} and
{\sc V.V. Peller}, {\em Functions of noncommuting self-adjoint operators under perturbation and estimates of triple operator integrals}, Advances in Mathematics {\bf295} (2016) 1--52.


\bibitem{APLa}  {\sc A.B. Aleksandrov} and {\sc V.V. Peller},
{\em Operator H\"older--Zygmund functions}, Advances in Math.
{\bf224} (2010), 910--966.

\bibitem{APSp}  {\sc A.B. Aleksandrov} and {\sc V.V. Peller},  {\em Functions of operators under perturbations of
class $\bS_p$}, J. Funct. Anal. {\bf258} (2010), 3675--3724.

\bibitem{APUnb}  {\sc A.B. Aleksandrov} and {\sc V.V. Peller}, {\em Functions of perturbed unbounded
self-adjoint operators. Operator Bernstein type inequalities}, Indiana Univ. Math. J. {\bf 59}
(2010), 1451--1490.

\bibitem{APJST}  {\sc A.B. Aleksandrov} and {\sc V.V. Peller}, {\em Trace formulae for perturbations of class $\boldsymbol{S}_m$}, J. Spectral Theory, {\bf1} (2011), 1--26.

\bibitem{APd} {\sc A.B. Aleksandrov} and {\sc V.V. Peller},
{\em Functions of perturbed dissipative operators}, Algebra i Analiz
{\bf23}:2 (2011), 9--51 (Russian); English transl., St. Petersburg Math. J., 
{\bf23}:2 (2012), 209--238.

\bibitem{APKomu}  {\sc A.B. Aleksandrov} and {\sc V.V. Peller},
{\em Functions of almost commuting operators and an extension of the Helton--Howe trace formula}.  J. Funct. Anal. {\bf271} (2016), 3300--3322.

\bibitem{APu} {\sc A.B. Aleksandrov} and {\sc V.V. Peller}, 
{\it Krein's trace formula for unitary operators and operator Lipschitz functions}  
Funkts. Anal. i Pril., {\bf50:3} (2016), 1--11 (Russian); English transl.,
Funct. Anal and Appl., {\bf50:3} (2016), 167--175.


\bibitem{APOL} {\sc A.B. Aleksandrov} and {\sc V.V. Peller}, {\it Operator Lipschitz functions}, 
Uspekhi Mat. Nauk {\bf71:4} (2016), 3--106; English transl.,
Russian Mathematical Surveys {\bf71:4} (2016), 605--702.

\bibitem{APMOI}  {\sc A.B. Aleksandrov} and {\sc V.V. Peller}, {\em Multiple operator integrals, Haagerup and Haagerup-like tensor products, and operator ideals}, Bull. London Math. Soc. {\bf49} (2017) 
463--479.


\bibitem{APCD}{\sc A.B. Aleksandrov} and {\sc V.V. Peller}, {\em Functions of perturbed commuting dissipative operators}, Math. Nachr., {\bf295} (2022), 1042--1062.

\bibitem{APdi} {\sc A.B. Aleksandrov} and {\sc V.V. Peller}, {\it Functions of perturbed pairs of noncommutative dissipative operators}, Algebra i Analiz, {\bf34:3} (2022), 93--114 (Russian); English transl., 
St. Petersburg Math. J. 34 (2023), 379--392.

\bibitem{APDok} {\sc A.B. Aleksandrov} and {\sc V.V. Peller}, {\it Functions of pairs of unbounded noncommuting self-adjoint operators under perturbation}, Dokl. Akad. Nauk, {\bf507} (2022), 5--9 (Russian); English transl.,
Dokl. Math. {\bf106:3} (2022), 407--411.

\bibitem{APne} {\sc A.B. Aleksandrov} and {\sc V.V. Peller}, {\it Functions of perturbed noncommuting unbounded self-adjoint operators}, Algebra i Analiz, {\bf34:6} (2022), 34--54 (Russian); English transl., 
St. Petersburg Math. J. 

\bibitem{APPS} {\sc A.B. Aleksandrov, V.V. Peller, D. Potapov}, and
{\sc F. Sukochev}, {\em Functions of normal operators under perturbations},
Advances in Math. {\bf226} (2011), 5216--5251.

\bibitem{ABB} {\sc J.M. Anderson, K.F. Barth} and {\sc D.A. Brannan}, {\em Research problems in complex analysis}, Bull. London Math. Soc. {\bf9} (1977), 129--162.

\bibitem{ACDS} {\sc N.A. Azamov, A.L. Carey, P.G. Dodds,} and {\sc F.A. Sukochev},
{\it Operator integrals, spectral shift and spectral flow}, 
Canad. J. Math. {\bf61} (2009), 241--263.

\bibitem{BGT}{\sc Ch. Batty, A. Gomilko, Yu. Tomilov}, {\em The theory of Besov functional calculus: developments and applications to semigroups}, J. Funct. Anal. 281 (2021), 109089.

\bibitem{Be} {\sc G. Bennett}, {\em Schur multipliers}, Duke Math. J. {\bf44} (1977), 603--639.

\bibitem{Bern} {\sc S.N. Bernstein}, {\em  A generalization of an inequality of S. B. Stechkin to entire functions of finite degree},
Doklady Akad. Nauk SSSR,  {\bf60}, 1948,  1487--1490 (Russian).

\bibitem{Bes}
{\sc O.V. Besov}, {\em On some families of functional spaces. Imbedding and extension theorems},  {\bf126:6} (1959), 1163--1165 (Russian).

\bibitem{BIN}
{\sc O.V. Besov, V.P. Il'in} and {\sc S.M. Nikolskii},  {\em Integral representations of functions, and embedding theorems}, Moscow, Nauka, 1996 (Russian).

\bibitem{BS1} {\sc M.S. Birman} and {\sc M.Z. Solomyak}, 
{\em Double Stieltjes operator integrals}, 
Problems of Math. Phys.,
Leningrad. Univ. 1 (1966), 33--67 (Russian). English transl., Topics Math. Physics 1, 25--54 . Consultants
Bureau Plenum Publishing Corporation, New York, 1967.

\bibitem{BS2} {\sc M.S. Birman} and {\sc M.Z. Solomyak}, 
{\em Double Stieltjes operator integrals. II}, Problems of Math. Phys.,
Leningrad. Univ. 2 (1967), 26--60 (Russian). English transl., Topics Math. Physics 2, 19--46. Consultants
Bureau Plenum Publishing Corporation, New York, 1968.

\bibitem{BS3} {\sc M.S. Birman} and {\sc M.Z. Solomyak}, 
{\em Double Stieltjes operator integrals. III},
Problems of Math. Phys. {\bf6}, Leningrad. Univ, 27--53 (1973) (Russian).

\bibitem{BS} {\sc M.S. Birman} and {\sc M.Z. Solomyak},  Spectral theory of selfadjoint operators in Hilbert space, Mathematics and its Applications (Soviet Series), D. Reidel Publishing Co., Dordrecht, 1987.

\bibitem{CP} {\sc R.W. Carey} and  {\sc J.D. Pincus},
{\em Mosaics, principal functions and mean motion in von Neumann algebras},
Acta Math. {\bf 138} (1977), 153--218. 

\bibitem{DK} {\sc Yu.L, Daletskii} and {\sc S.G. Krein}, 
{\it Integration and differentiation of functions of Hermitian operators and
application to the theory of perturbations}, Trudy Sem. Functsion. Anal., Voronezh. Gos.
Univ. 1 (1956), 81--106  (Russian).

\bibitem{EL} {\sc P.P.B. Eggermont} and {\sc Y.J. Leung}, {\em On a factorization problem for convergent sequences and on Hankel forms in bounded sequences}, Proc. Amer. Math. Soc. {\bf96} (1986), 269--274.

\bibitem{FaO} {\sc Yu. B. Farforovskaya}, {\em The connection of the Kantorovich-Rubinstein metric for spectral resolutions of selfadjoint operators with functions of operators},
Vestnik Leningrad. Univ., {\bf 19} (1968), 94--97 (Russian).


\bibitem{Gar} {\sc J. Garnett}, {\em Bounded analytic functions}, Academic Press, Inc., New York-London, 1981.

\bibitem{GK} {\sc I.C. Gohberg} and {\sc M.G. Krein}, {\em Introduction to the theory of linear nonselfadjoint operators}, Nauka, Moscow, 1965 (Russian); English transl., Translations of Mathematical Monographs, {\bf18}, American Mathematical Society, Providence, RI, 1969.

\bibitem{Gr} {\sc A. Grothendieck}, {\em R\'esum\'e de la th\'eorie m\'etrique des produits tensoriels
topologiques}, Bol. Soc. Matem. Sao Paulo {\bf8} (1955), 1--79.

\bibitem{Ha} {\sc M. Haase}, {\em Transference principles for semigroups and a theorem of Peller}, J. Funct. Anal. {\bf261} (2011), 2959--2998.

\bibitem{HH} {\sc J.W. Helton} and {\sc R. Howe},
{\em Integral operators, commutators, traces, index, and homology}, 
in ``Lecture Notes in Math.'', vol. {\bf345}, pp. 141--209,
Springer-Verlag, New York, 1973.

\bibitem{IS1} {\sc I.A. Ibragimov} and {\sc V.N. Solev}, {\em Certain analytic problems arising in the theory of stationary random processes}, In: "99 open problems in linear and complex analysis", Zapiski. Nauchn. Semin. LOMI {\bf81} (1978), 70--72 (Russian).

\bibitem{IS2} {\sc I.A. Ibragimov} and {\sc V.N. Solev}, {\em A certain condition for the regularity of Gaussian stationary sequence}, Zapiski. Nauchn. Semin. LOMI  {\bf12} (1969), 113--125 (Russian).

\bibitem{Ko} {\sc L.S. Koplienko}, {\em The trace formula for perturbations of nonnuclear type},
Sibirsk. Mat. Zh., {\bf25:5} (1984),  62--71 (Russian); English transl., Siberian Math. J. {\bf25} (1984), 735--743.

\bibitem{JTT} {\sc K. Juschenko, I.G. Todorov} and {\sc L. Turowska}, {\em Multidimensional operator multipliers}, Trans. Amer. Math. Soc. {\bf361}
(2009), 4683--4720.






\bibitem{KS1} {\sc E. Kissin} and {V.S. Shulman}, 
{\em On a problem of
J. P. Williams}. Proc. Amer. Math. Soc. {\bf130} (2002), 3605 -- 3608.

\bibitem{KS2} {\sc E. Kissin} and {V.S. Shulman},
{\em Classes of operator-smooth functions. I. Operator-Lipschitz functions},
Proc. Edinb. Math. Soc. (2) {\bf48} (2005), 151--173.

\bibitem{KS3} {\sc E. Kissin} and {V.S. Shulman},
{\em Classes of operator-smooth functions. II. Operator-differentiable functions},
Int. Equat. Oper. Theory {\bf49} (2004), 165--210.

\bibitem{KS4} {\sc E. Kissin} and {V.S. Shulman},
{\em On fully operator Lipschitz functions}.
J. Funct. Anal. {\bf253} (2007), 711--728.


\bibitem{Kr1} {\sc M.G. Krein}, {\em On the trace formula in perturbation theory},
Mat. Sbornik N.S. {\bf33} (1953), 597--626 (Russian).

\bibitem{Kr2} {\sc M.G. Krein}, {\em On perturbation determinants and a trace formula for unitary and self-adjoint operators}, Dokl. Akad. Nauk SSSR, {\bf144:2} (1962), 268--271 (Russian).

\bibitem{Kr} {\sc M.G. Krein}, {\em Certain new Banach algebras and theorems of the type of the Wiener-L?vy theorems for series and Fourier integrals}, Mat. Issled., Math. Inst. Acad. Nauk Mold. SSR {\bf1:1} (1966), 82--109 (Russian).

\bibitem{Kr87}
{\sc M.G. Krein}, {\em Perturbation determinants and a trace formula for some classes of pairs of operators},
J. Operator Th., {\bf17} (1987), 129--187.

\bibitem{Ku} {\sc P. Koosis}, {\em Introduction to $H_p$ spaces. With an appendix on Wolff's proof of the corona theorem}, London Mathematical Society Lecture Note Series, 40. Cambridge University Press, Cambridge-New York, 1980.

\bibitem{Kro} {\sc L. Kronecker}, {\em Zur Theorie der Elimination einer Variablen aus zwei algebraischen Gleichungen}, Monastber. K\"onigl. Preuss. Akad. Wiss., Berlin, 1981, 535--600.

\bibitem{La} {\sc H. Langer}, {\it Eine Erweiterung der Spurformel 
der St\"orungstheorie}, Math. Nachr. {\bf30} (1965), 123--135.

\bibitem{Leb} {\sc A. Lebow}, {\it A power bounded operator that is not polynomially bounded}, Michigan Math. J. {\bf15} (1968), 397--399.

\bibitem{Li}{\sc I.M. Lifshits}, {\em On a problem of the theory of perturbations connected with quantum statistics}, Uspekhi Mat. Nauk {\bf 7:1(47)} (1952), 171--180 (Russian).

\bibitem{KwaP}{\sc S. Kwapie\'n} and {\sc A. Pe\l czy\'nski}, {\em On two problems of S. Mazur from the Scottish book}, lecture at the Colloquium dedicated to the memory of Stanis\l aw Mazur, Warsaw Univ. 1985 (unpublished).

\bibitem{Lev}
{\sc B.~Ya. Levin}, {\em Lectures on entire functions},
Translation of Math. Monogr., AMS, vol. 150, 1996.

\bibitem{MN}
{\sc M. Malamud} and {\sc H. Neidhardt},
{\it Trace formulas for additive and non-additive perturbations},
Adv. Math. {\bf274} (2015), 736--832.

\bibitem{MNP1} 
{\sc M. M. Malamud, H. Neidhardt} and {\sc V. V. Peller}, {\it Analytic operator Lipschitz functions in the disk and a trace formula for functions of contractions}, Funktsional. Anal. i Prilozhen., {\bf51:3} (2017), 33--55
(Russian); English transl., Funct. Anal. Appl. {\bf51} (2017), 185--203.

\bibitem{MNP2}
{ \sc M.M. Malamud, H. Neidhardt} and {\sc V.V. Peller}, {\it Absolute continuity of spectral shift}, J. Funct. Anal., {\bf276} (2019), 1575--1621.


\bibitem{Ma} {\sc R.D. Mauldin}, {\em The Scottish Book}, Birkh\"auser, 2015.

\bibitem{McS} {\sc E. McDonald} and {\sc F. Sukochev},
{\it Lipschitz estimates in quasi-Banach Schatten ideals}, 
Math. Ann. {\bf383} (2022), 571--619. 

\bibitem{NP} {\sc F.L. Nazarov} and {V.V. Peller}, {\em Functions of $n$-tuples of commuting self-adjoint operators}, J. Funct. Anal. {\bf266} (2014), 5398--5428.

\bibitem{Nei} {\sc H. Neidhardt}, {\em Spectral shift function and Hilbert--Schmidt perturbation: extensions of some work of L.S. Koplienko}, Math. Nachr. {\bf138} (1988), 7--25.

\bibitem{NF}  {\sc L. N.  Nikol'skaya} and {\sc Yu. B. Farforovskaya},
{\em Operator H\"oldericity of H\"older functions}, Algebra i Analiz {\bf22:4}
(2010), 198 -- 213 (Russian); English transl., St. Petersburg Math. J. {\bf22} (2011), 657--668 


\bibitem{Nik} {\sc N.K. Nikolski}, {\em Operators, functions, and systems: an easy reading. Vol. 1. Hardy, Hankel, and Toeplitz}, Mathematical Surveys and Monographs, {\bf92}. 
American Mathematical Society, Providence, RI, 2002.

\bibitem{SNi} {\sc S.M. Nikol'skii}, {\it Approximation of functions of several variables and imbedding theorems}, Moscow, Nauka, 1969 (Russian); English transl.,  Die Grundlehren der mathematischen Wissenschaften, Band 205. Springer-Verlag, New York-Heidelberg, 1975.

\bibitem{Pee} {\sc J. Peetre},
{\em New thoughts on Besov spaces}, Duke Univ. Press., Durham, NC, 1976.

\bibitem{Pek} {\sc A. A. Pekarskii},
{\em Classes of analytic functions defined by best rational approximations in $H_p$}, 
Mat. Sb. (N.S.) {\bf127(169)} (1985), no. 1, 3--20, 142
(Russian); English transl., Math. USSR-Sb. {\bf55} (1986), 1--18.

\bibitem{PS} {\sc A. Pe\l czy\'nski} and {\sc F. Sukochev}, {\it Some remarks on Toeplitz multipliers and Hankel matrices}, Studia Math. {\bf175} (2006), 175--204.

\bibitem{PeSb1} {\sc V.V. Peller}, {\em Hankel operators of class ${\frak S}_p$ and their applications (rational approximation, Gaussian processes, the problem of majorization of operators).},  Mat. Sb. (N.S.),
{\bf 113(155):4(12)} (1980), 538--581 (Russian); English transl., Math. USSR-Sb., {\bf41:4} (1982), 443--479.

\bibitem{Pe82} {\sc V.V. Peller}, {\ Estimates of functions of power bounded operators on Hilbert spaces}, J. Operator Theory {\bf7} (1982), 341--372.

\bibitem{PeSb2} {\sc V.V. Peller}, {\em Description of Hankel operators of the class ${\frak S}_p$ for $p>0$, investigation of the rate of rational approximation and other applications},  Mat. Sb. (N.S.)
{\bf 122(164):4(12)} (1983), 481--510 (Russian); English transl., Math. USSR-Sb., {\bf50}:2 (1985), 465--494.

\bibitem{PeFAP} {\sc V.V. Peller}, {\em Hankel operators in the theory of perturbations of unitary and selfadjoint operators},  Funktsional. Anal. i Prilozhen. {\bf19:2}  (1985), 37--51 (Russian); English transl., Funct. Anal. Appl., {\bf19:2} (1985), 111--123.

\bibitem{PeSzh} {\sc V.V. Peller},
{\em For which $f$ does $A-B\in{\bf S}_{p}$
imply that $f(A)-f(B)\in{\bf S}_{p}$?}, Operator Theory, Birkh\"{a}user,
{\bf 24} (1987), 289--294.

\bibitem{PeWH} {\sc V.V. Peller}, {\em Wiener--Hopf operators on a finite interval and Schatten--von Neumann classes}, Proc. Amer. Math. Soc. {\bf104} (1988), 479--486.

\bibitem{PeMD90} {\sc V.V. Peller} {\em Hankel operators in the perturbation 
theory of unbounded self-adjoint operators}.  Analysis and partial differential equations,  529--544,
Lecture Notes in Pure and Appl. Math., {\bf122}, Dekker, New York, 1990.

\bibitem{PeAC} {\sc V.V. Peller},
{\em Functional calculus for a pair of almost
commuting selfadjoint operators}, J. Funct. Anal., {\bf 112} (1993), 325--345.

\bibitem{Pe} {\sc V. V. Peller},  
{\em Hankel operators and their applications},
Springer Monogr. Math., Springer-Verlag, New York,
2003.

\bibitem{PeKN} {\sc V.V. Peller}, {\em An extension of the Koplienko--Neidhardt trace formulae}, 
 J. Funct. Anal. {\bf221} (2005), 456--481.
 
 \bibitem{PeOD} {\sc V.V. Peller}, {\em Multiple operator integrals and higher operator
derivatives}, J. Funct. Anal.  {\bf233}  (2006),  515--544.

\bibitem{PeStud}  {\sc V.V. Peller}, {\em On S. Mazur's problems 8 and 88 from the Scottish Book}, 
Studia Math. {\bf180 (2)} (2007), 191--198.

\bibitem{PeDifco}  {\sc V.V. Peller}, {\em Differentiability of functions of contractions}, 
Amer. Math. Soc, Transl. (2) {\bf226} (2009), 109--131.

\bibitem{PeSpSd} {\sc V.V. Peller}, {\it The Lifshits--Krein trace formula and operator Lipschitz functions}, Proc. Amer. Math. Soc. {\bf144} (2016), 5207--5215.

\bibitem{PeBMS} {\sc V.V. Peller}, {\em Multiple operator integrals in perturbation theory}, Bull. Math. Sci. {\bf6} (2016), 15--88.

\bibitem{Petri} {\sc V.V. Peller}, {\it Functions of triples of noncommuting self-adjoint operators under perturbations of class $\boldsymbol{S_p}$}, 
Proc. Amer. Math. Soc. {\bf146} (2018), 1699--1711.

\bibitem{PeCS} {\sc V.V. Peller}, {\em Functions of commuting contractions under perturbation},
Math. Nachr. {\bf 292} (2019), 1151--1160.

\bibitem{PeKh} {\sc V.V. Peller} and {\sc S.V. Khrushch\"ev}, {\em Hankel operators, best approximations and stationary Gaussian processes}, Uspekhi Mat. Nauk {\bf37:1(223)} (1982), 53--124 (Russian); English transl.,
Russian Math. Surveys 37 (1982), 61--144.

\bibitem{Pin1} {\sc J.D. Pincus}, {\em Commutators and systems of singular integral equations, I}, Acta Math. {\bf121} (1968), 219-249.

\bibitem{Pin2} {\sc J.D. Pincus}, {\em On the trace of commutators in the algebra of operators generated by an operator with trace class self-commutator}, Stony Brook preprint (1972).

\bibitem{PSS} {\sc D. Potapov, A. Skripka}, and {F. Sukochev},
{\em Spectral shift function of higher order}, Invent. Math. {\bf193} (2013), 501--538.

\bibitem{PSS2} {\sc D. Potapov, A. Skripka}, and {F. Sukochev}, {\it Higher-order spectral shift for contractions}, Proc. Lond. Math. Soc. {\bf(3) 108} (2014), 327--349.
 
\bibitem{PSS3} {\sc D. Potapov, A. Skripka}, and {F. Sukochev},
{\it Functions of unitary operators: derivatives and trace formulas}, J. Funct. Anal. {\bf270} (2016), no. 6, 2048--2072.

\bibitem{PoSu} {\sc D. Potapov} and {\sc F. Sukochev}, {\em Operator-Lipschitz functions in Schatten--von Neumann classes}, Acta Math. {\bf207} (2011), 375--389.

\bibitem{PoSu2} {\sc D. Potapov} and {\sc F. Sukochev}, {\it Koplienko spectral shift function on the unit circle}, Comm. Math. Phys. {\bf309} (2012) 693--702.


\bibitem{RSN} {\sc F. Riesz} and {\sc B. Sz.-Nagy}, {\em Le\c cons d'analyse fonctionnelle}, Gauthier-Villars, Paris; Akad'emiai Kiad\'o, Budapest, 1955.

\bibitem{Roc} {\sc R. Rochberg}, {\em Toeplitz and Hankel operators on the Paley--Wiener space}, Int. Equat. Oper. Theory {\bf10} (1987), 187--235.

\bibitem{Ro} {\sc M. Rosenblum}, {\em The absolute continuity of Toeplitz's matrices}, Pacific J. Math. {\bf10} (1960), 987--996.

\bibitem{Ryb1} {\sc A. V. Rybkin}, {\em The spectral shift function for a dissipative and a selfadjoint operator, and trace formulas for resonances}, Mat. Sb. (N.S.), {\bf125:3} (1984),  420--430. 

\bibitem{Ryb2} {\sc A. V. Rybkin}, {\em A trace formula for a contractive and a unitary operator},  Funktsional. Anal. i Prilozhen., {\bf21:4} (1987), 85--87 (Russian). 

\bibitem{Ryb3} {\sc A. V. Rybkin}, {\em The discrete and the singular spectrum in the trace formula for a contractive and a unitary operator.},  Funktsional. Anal. i Prilozhen., {\bf23:3} (1989), 84--85 (Russian); English transl., Funct. Anal. Appl. {\bf23:3} (1989), 244--246.

\bibitem{Ryb4} {\sc A. V. Rybkin}, {\em The spectral shift function, the characteristic function of a contraction and a generalized integral.}, Mat. Sb., {\bf185:10} (1994),  91--144 (Russian); English transl.,  Russian Acad. Sci. Sb. Math. {\bf83} (1995), 237--281.

\bibitem{Sch} {\sc L. Schwartz}, {\em Th\'eorie des distributions}, Paris, Hermann, 1966.

\bibitem{Se} {\sc S. Semmes}, {\em Trace ideal criteria for Hankel operators and applications to Besov spaces}, Int. Equat. Oper. Theory {\bf7} (1984), 241--281.

\bibitem{SNF} {\sc B. Sz.-Nagy} and {\sc C. Foia\c s}, {\em Analyse harmonique des op\'erateurs de l'espace de Hilbert}, Masson et Cie, Paris; Akad\'emiai Kiad\'o, Budapest, 1967.


%


















\bibitem{T} {\sc H. Triebel}, {\em Theory of function spaces}, Birkh\"auser Verlag, Basel, 1983.

\bibitem{V}{\sc P. Vitse}, {\em A Besov class functional calculus for bounded holomorphic semigroups}, 
J. Funct. Anal. {\bf228} (2005), 245--269.

\bibitem{Wi} {\sc J.P. Williams},  
{\em Derivation ranges: open problems},
In: Topics in modern operator theory, 319--328, Operator Theory: Adv. Appl., {\bf2}, Birkh\"auser, Basel-Boston, Mass., 1981. 


\end{thebibliography}
